\documentclass[11pt,a4paper]{article}
\pdfoutput=1
\usepackage[utf8]{inputenc}
\usepackage[T1]{fontenc}
\usepackage{amsmath,amssymb,amsthm}
\usepackage{mathtools}
\usepackage{geometry}
\usepackage{hyperref}\hypersetup{pdftitle={Shift-Symmetric Configurations in d-Dimensional Cellular Automata: Enumeration, Bounds, Irreversibility, and Detection},pdfauthor={Peter Banda, John Caughman, Tobi Lehman}}
\usepackage{booktabs}
\usepackage{array}
\usepackage{enumitem}
\usepackage{graphicx}
\usepackage{tikz}
\usetikzlibrary{positioning,arrows.meta,calc}

\newtheorem{theorem}{Theorem}[section]
\newtheorem{lemma}[theorem]{Lemma}
\newtheorem{corollary}[theorem]{Corollary}
\newtheorem{proposition}[theorem]{Proposition}
\newtheorem{definition}[theorem]{Definition}
\theoremstyle{remark}
\newtheorem{remark}[theorem]{Remark}
\newtheorem{observation}[theorem]{Computational Observation}

\DeclareMathOperator{\lcm}{lcm}
\DeclareMathOperator{\Fix}{Fix}
\DeclareMathOperator{\Stab}{Stab}
\DeclareMathOperator{\rk}{rk}

\newcommand{\ZZ}{\mathbb{Z}}

\newcommand{\FF}{\mathbb{F}}
\newcommand{\A}{\mathcal{A}}
\newcommand{\vv}{\mathbf{v}}
\newcommand{\uu}{\mathbf{u}}
\newcommand{\ww}{\mathbf{w}}
\newcommand{\ssb}{\mathbf{s}}
\newcommand{\gaussbinom}[2]{\genfrac{[}{]}{0pt}{}{#1}{#2}}

\title{Shift-Symmetric Configurations in $d$-Dimensional Cellular Automata:\\
Enumeration, Bounds, Irreversibility, and Detection}

\author{Peter Banda\thanks{Author to whom correspondence should be addressed: \texttt{peter.banda@protonmail.com}}\\{\normalsize Cequence} \and John Caughman\\{\normalsize Portland State University} \and Tobi Lehman\\{\normalsize Functorial, Inc.}}
\date{}

\begin{document}
\maketitle

\begin{abstract}
A configuration of a cellular automaton on the $d$-dimensional torus
$\ZZ_n^d$ is \emph{shift-symmetric} if it is fixed by some nonzero
translation. Banda, Caughman, Cenek, and Teuscher (\emph{Chaos} 29, 063120,
2019) enumerated these configurations in two dimensions and drew dynamical
and algorithmic consequences. We extend their results to every dimension.
We give an explicit closed-form enumeration for all $n$, $d$, and alphabet
sizes $\sigma$, in which the two-dimensional coefficients $1$, $p{+}1$, and
$p$ become $r_d(p,v)=p^{\binom v2}\gaussbinom dv_p$; it extends to
configurations with a prescribed number of active cells and to all finite
abelian groups. The case of prime
$n$ recovers, in complement form, a result of Castillo-Ramirez and
Gadouleau; the multi-prime closed form and the dynamical consequences
beyond the stabilizer inclusion are new. We prove a lower bound attained
\emph{exactly} when $n$ is prime, and upper bounds, among them
$(2^d-1+2\log_2 n)\,\sigma^{n^d/2}$. As $d$ grows, the symmetry probability
decays as $\sigma^{-c\,n^d}$ with $c\to1-1/p_1(n)$, set solely by the
smallest prime factor of $n$, and the symmetry advantage of the side
$p{+}1$ over a prime side $p$ inverts within an explicit window around
$d_c(p)\approx p\ln2$. Under the necessary hypothesis of
translation-equivariant neighborhoods, shift-symmetry is irreversible, the
per-prime symmetry ranks never decrease, and non-symmetric tasks are
unsolvable from symmetric inputs. We give a detection algorithm with
worst-case cost $O\bigl(n^d\sum_{p\mid n}\frac{p^d-1}{p-1}\bigr)$ and an
average-case analysis, exact for each subgroup test. The mathematical
content of every numbered statement is formally verified in Lean~4, and the
correspondence between each theorem, lemma, proposition, and corollary and
its formal statement was audited by hand.
\end{abstract}
%
%

\begin{quote}
\bfseries
Cellular automata are grids of cells that update in parallel, every cell
applying the same local rule to its neighbors at the same relative
positions. On a grid that wraps around like a torus, some patterns are
unchanged when shifted by a nonzero offset. Such a uniform update can never
destroy this symmetry, so tasks that require breaking it, such as electing
a single leader cell, cannot be solved from these patterns. How likely is a
random starting pattern to be one of them, and how does this depend on the
grid? For two-dimensional grids, these questions were answered in 2019.
Here we give an exact formula in every dimension. The arithmetic of the
side length governs the answer: a simple formula bounds the count from
below and is exact precisely when the side length is prime, and as the
dimension grows, the smallest prime factor of the side length alone sets
the rate per cell at which such patterns become rare. All theorems were
also verified by computer in the Lean proof assistant.
\end{quote}

\newpage
\tableofcontents
\newpage

\section{Introduction}

\begin{figure}[t]
\centering
\includegraphics[width=\textwidth]{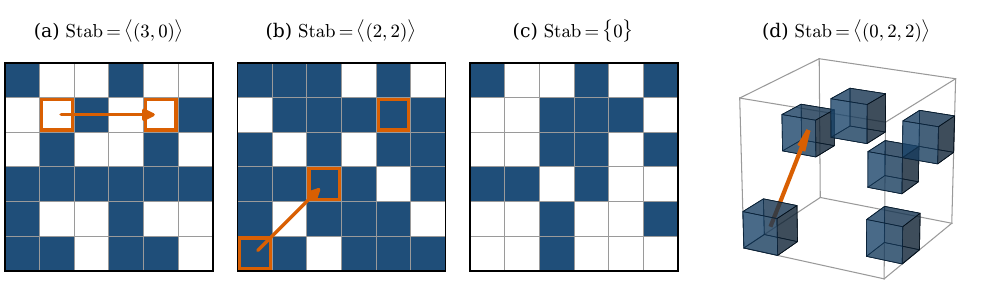}
\caption{Shift-symmetric and asymmetric configurations ($\sigma=2$: filled
and empty cells). A cell $\uu=(u_1,u_2)$ is drawn in column $u_1$ and row
$u_2$, counted from the bottom left, and $\ssb$ is fixed by $\vv$ when
$\ssb(\uu+\vv)=\ssb(\uu)$ for all $\uu$; each arrow runs from a cell $\uu$
to $\uu+\vv$ for a generator $\vv$ of the stabilizer. (a)~A configuration on
$\ZZ_6^2$ with stabilizer $\langle(3,0)\rangle$: the translation by $(3,0)$
maps every cell to a cell in the same state (outlined: one orbit).
(b)~Stabilizer $\langle(2,2)\rangle$, a diagonal translation of order~$3$.
(c)~Trivial stabilizer: this configuration is not shift-symmetric. (d)~A
configuration on $\ZZ_4^3$, its six filled cells drawn, with stabilizer
$\langle(0,2,2)\rangle$; the arrow runs from $(0,0,0)$ to $(0,2,2)$. Of the
$2^{36}$ configurations on $\ZZ_6^2$, exactly $801\,064$ (about
$1.2\cdot10^{-5}$ of them) are shift-symmetric (Theorem~\ref{thm:main-enum};
a reference value of Section~\ref{sec:compver}).}
\label{fig:examples}
\end{figure}

A configuration of a cellular automaton on the $d$-dimensional torus
$\ZZ_n^d$ is said to be \emph{shift-symmetric} if it is fixed by some nonzero
translation (Figure~\ref{fig:examples}).
For the two-dimensional case, Banda, Caughman, Cenek, and Teuscher
(\emph{Chaos} 29, 063120, 2019) studied shift-symmetry for configurations of
two-dimensional cellular automata (CA) on toroidal lattices
$\ZZ_n\times\ZZ_n$. Such symmetry is noteworthy since,  once a configuration is fixed by a nonzero translation, no uniform synchronous update can ever cause the symmetry to be lost. As a result, symmetric
configurations effectively function as irreversible traps that can make non-symmetric tasks
(leader election, pattern recognition) unsolvable in the worst case.
Quantifying this obstruction requires counting the shift-symmetric
configurations, which Banda et al.~\cite{banda2019} accomplished for $d=2$, leaving higher-dimensional extensions for
future consideration. In this paper, we provide the $d$-dimensional extension
and, in the process, sharpen and correct several statements whose naive
generalizations to higher dimensions fail.

\paragraph{Contributions.} The primary contributions of this paper are the following:
\begin{enumerate}[label=(\roman*)]
    \item An explicit closed-form enumeration of the set $S_{n^d}$ of shift-symmetric
    configurations on $\ZZ_n^d$ for arbitrary $n$, $d$, and alphabet size
    $\sigma$ (Theorem~\ref{thm:main-enum}), with two independent
    derivations, and a refinement to configurations with exactly $k$
    active cells (Theorem~\ref{thm:kactive}); both extend to every finite
    abelian group, in particular to rectangular tori, with $d$ replaced by
    the $p$-ranks (Theorem~\ref{thm:abelian}, Remark~\ref{rem:abelian}), and
    the count is smallest for
    the cyclic group among the abelian groups of a given order
    (Corollary~\ref{cor:abelian-ranks}).
    \item A lower bound for $|S_{n^d}|$  attained exactly when $n$ is prime
    (Theorem~\ref{thm:lower-bound}). We show that a naive generalization of the
    two-dimensional bound of \cite{banda2019} is \emph{false} for $d\ge3$
    (Remark~\ref{rem:oldbound-false}). We also give an upper bound with an
    explicit prefactor $2^d-1+2\log_2 n$ (Theorem~\ref{thm:upper-bound}).
    \item Proofs of irreversibility, insolvability, and per-prime symmetry-rank
    monotonicity in $d$ dimensions (Section~\ref{sec:dynamics}), under a
    translation-equivariance hypothesis that we show is \emph{necessary}: a
    two-cell counterexample (Remark~\ref{rem:equivariance-needed})
    demonstrates that a shared local rule with cell-dependent neighborhoods
    does not preserve symmetry.
    \item A detection algorithm for shift-symmetry with worst-case bound
    $O(\mathrm{sopf}_d(n)\cdot n^d)$ and an exact average-case expectation
    analysis (Section~\ref{sec:detection}), with the average-case
    ``polylogarithmic'' claims of the two-dimensional literature restated
    under appropriate hypotheses.
    \item A number of computational observations regarding higher
    dimensions (Section~\ref{sec:empirical}): the prime/composite extremal
    pattern of $d=2$ \emph{breaks} for $d\ge3$, and the multi-term
    correction regime collapses rapidly with growing $p$ and $d$.
\end{enumerate}

\paragraph{Comments on verification.} Every enumerative formula below has been checked, with exact integer
arithmetic, against definition-level exhaustive enumeration and against an
independent subgroup-level inclusion--exclusion evaluator; the two main
theorems additionally passed independent adversarial review. The
supplementary material, namely the scripts that reproduce every
enumeration, scan, and figure (with their domains, outputs, and
checksums),\footnote{Catalogued in \texttt{research/experiments/}; see the
reproducibility paragraph in Section~\ref{sec:empirical}. A companion
erratum to \cite{banda2019}, correcting misprints and gaps in its proofs, is
in \texttt{research/2d\_fix\_note/}.} the Lean
development (Section~\ref{sec:formal}), and the correspondence audit
(Section~\ref{sec:conclusion}), is in the project repository,
\href{https://github.com/peterbanda/d-dim-shift-symmetry-in-ca}{\nolinkurl{github.com/peterbanda/d-dim-shift-symmetry-in-ca}}.
Beyond testing and
review, the numbered results are \emph{machine-checked}:
Theorems~\ref{thm:main-enum} and~\ref{thm:lower-bound}, in the generality
stated here, and the mathematical content of the other numbered
statements of Sections~\ref{sec:model}--\ref{sec:empirical} and
Appendices~\ref{app:lower} and~\ref{app:abelian} (statistical summaries aside) have been formally
verified in the Lean~4 proof assistant (Section~\ref{sec:formal}).

\paragraph{Related work.}
For the two-dimensional torus, Banda et al.~\cite{banda2019} enumerated shift-symmetric
configurations, bounded their number, and analyzed detection; to our
knowledge no extension has appeared since. Earlier, CA leader election had
been studied in one dimension, both by evolving binary rules \cite{banda2011}
and through symmetry-based impossibility results and performance upper
bounds \cite{banda2015}. From the group-theoretic side,
Castillo-Ramirez and Gadouleau \cite{cg2017} counted configurations with a
given stabilizer class over an \emph{arbitrary} finite group via the
M\"obius function of the subgroup lattice. Their Lemma~7, which uses Hall's
evaluation of the M\"obius function of a $p$-group (\cite{hall1933}; see
Hawkes, Isaacs, and \"Ozaydin \cite[Cor.~3.5]{hio1989}; see also Delsarte
\cite{delsarte1948} for M\"obius inversion over the subgroups of finite
abelian groups), yields in
complement form the \emph{prime-power} case $n=p^\alpha$ of our enumeration
theorem: by its part~(iii), it remains to count the elementary abelian
subgroups of $\ZZ_{p^\alpha}^d$, which are the subspaces of its $p$-torsion
$\cong\FF_p^d$ (Lemma~\ref{lem:torsion}); part~(iv) is the prime case. Section~\ref{sec:second-derivation}
shows how the general case can alternatively be derived from their
framework. What appears neither there
nor in the least-period counts of Castillo-Ramirez and S\'anchez-\'Alvarez
\cite{cs2022}, whose worked examples stop at rank two, is the explicit
multi-prime closed form for $\ZZ_n^d$ and for finite abelian groups
(Section~\ref{sec:abelian}), the prime-tightness characterization for
$d\ge2$ (for $d=1$ it is noted after \cite[Lemma~7]{cg2017}), or the
consequences drawn in
Section~\ref{sec:dynamics} (insolvability, rank monotonicity). On the
combinatorial side, Ethier \cite{ethier2013} and Ethier and Lee
\cite{ethierlee2015} count two-dimensional toroidal arrays up to symmetry
(orbits, not symmetric configurations), and Adamson, Deligkas, Gusev, and
Potapov study multidimensional necklaces algorithmically and state a
recursion across dimensions for the number of ``atranslational'' ones ---
the orbit-level complement of our shift-symmetric configurations --- on
rectangular grids \cite[Thm.~3]{adamson2021};\footnote{As printed, this
recursion does not reproduce exhaustive counts: none of the readings of its
formula that we tried matches them on all grids with at most $16$ cells over
two and three symbols (Run~037 of the repository).} they derive no non-recursive
formula for this count (it is $(\sigma^N-|S_G|)/N$ with $|S_G|$ from
Theorem~\ref{thm:abelian}; see Remark~\ref{rem:abelian}; see also the
group-theoretic necklace polynomials of Oh \cite{oh2012}) and no dynamical
consequences. In symbolic dynamics, the numbers of points fixed by a
\emph{given} finite-index subgroup $L\le\ZZ^d$ (for the full shift,
$\sigma^{[\ZZ^d:L]}$) are classical and are packaged by Lind's zeta function
for $\ZZ^d$-actions \cite{lind1996}; Miles and Ward \cite{mw2009} invert
them over the lattice of finite-index subgroups to count the orbits with a
given stabilizer, and study the growth of these counts, without reducing
them to the explicit rank-indexed formula used here. Our object is, instead, the union over all nontrivial
translations on a fixed finite torus. Lower
bounds on the number of aperiodic configurations of a general finite group
were given by Gao, Jackson, and Seward \cite[Cor.~1.7.2]{gjs2016} and
sharpened in \cite[Thm.~5]{cg2017}. On $\ZZ_n^d$, the former is
exponentially weaker in the exponent than Proposition~\ref{prop:rate};
the latter, whose proof contains our union bound
(Theorem~\ref{thm:upper-bound}(i)), already implies its limit, with a
larger prefactor for $n\ge7$. For
$d=1$ our formulas reduce to classical aperiodic-word counts going back to
Moreau and Witt \cite{moreau1872,witt1937}; see also
\cite{metropolisrota1983}. The prime-lattice, binary, two-dimensional case
of our count also appears, without proof, as a formula line in OEIS entry
A323863 \cite{oeisA323863} (orbit-level companion entries: A323872,
A184271). On the dynamical side,
Kornyak \cite{kornyak2010} observed that deterministic discrete dynamics
cannot decrease symmetry, and for cellular automata on a finite group the
inclusion of the stabilizer of a configuration in that of its image
(Theorem~\ref{thm:irrev}) is stated in the proof of \cite[Lemma~10]{cg2017},
with equality for invertible ones in that of \cite[Lemma~3]{cg2017}; its
converse (some cellular automaton maps $x$ to $y$ if $\Stab(x)\le\Stab(y)$)
is \cite[Lemma~10]{cg2017} for distinct orbits and \cite[Thm.~2(iii)]{cgvn2017}
in general. The orbit count of Remark~\ref{rem:abelian}(ii) is the case
$H=\{e\}$ of \cite[Lemma~4]{cg2017}.

\section{Model}\label{sec:model}

In this section we introduce some of the basic definitions and notation we will use in our study of cellular automata and their configurations.

\begin{definition}[Lattice and configurations]
Fix integers $n\ge2$, $d\ge1$ and a finite alphabet $\Sigma$ with
$\sigma=|\Sigma|\ge2$ (for counting statements, $\sigma=1$ is also allowed).
We define the {\bf lattice} $\A$ to be the finite abelian group $\A=(\ZZ/n\ZZ)^d$, and we refer to the $N=n^d$ elements of $\A$ as the {\bf cells} of $\A$.
A {\bf configuration} is any function $\ssb:\A\to\Sigma$. 
\end{definition}

\begin{definition}[Translations acting on configurations] For
$\vv\in\A$, the {\bf translation} $\tau_\vv$ acts on the set of configurations by
$(\tau_\vv\ssb)(\uu)=\ssb(\uu+\vv)$, for any $\uu, \vv \in \A$. For any configuration $\ssb$, the corresponding {\bf stabilizer}, denoted
$\Stab(\ssb)=\{\vv\in\A:\tau_\vv\ssb=\ssb\}$, is a subgroup of $\A$. For any subgroup $H\le\A$, its {\bf fixed set} is given by $\Fix(H)=\{\ssb:\tau_{\mathbf h}\ssb=\ssb\text{ for all }\mathbf h\in H\}$.
\end{definition}

\begin{definition}[Shift-symmetry] For ${\mathbf 0} \neq \vv \in \A$, we define
\[S_{n^d}(\vv)=\{\ssb:\tau_\vv\ssb=\ssb\}=\Fix(\langle\vv\rangle).\] 
The set of {\bf shift-symmetric configurations} is
\[
    S_{n^d} \;=\; \bigcup_{\mathbf 0\neq\vv\in\A} S_{n^d}(\vv)
    \;=\; \{\ssb : \Stab(\ssb)\neq\{\mathbf 0\}\}.
\]
\end{definition}

\begin{definition}[Equivariant cellular automaton]\label{def:ca}
A {\bf cellular automaton} on $\A$ over $\Sigma$ is given by a finite list
of {\bf offsets} $\mathbf o_1,\dots,\mathbf o_b\in\A$ and a {\bf local
rule} $\phi:\Sigma^b\to\Sigma$; the {\bf global map}
$\Phi:\Sigma^{\A}\to\Sigma^{\A}$ updates all cells synchronously:
\[
    \Phi(\ssb)(\uu)\;=\;\phi\bigl(\ssb(\uu+\mathbf o_1),\dots,
    \ssb(\uu+\mathbf o_b)\bigr).
\]
Equivalently, $\Phi$
is {\bf translation-equivariant}: $\Phi(\tau_\vv\ssb)=\tau_\vv\Phi(\ssb)$
for all $\vv\in\A$. 
\end{definition}

\begin{remark}[Equivariance is necessary, not cosmetic]
\label{rem:equivariance-needed}
We note that the property of translation-equivariance is used extensively in our discussion of dynamical results. Indeed, definitions in the style of ``a neighborhood function
$\eta:\A\times\Sigma^{\A}\to\Sigma^b$ together with a shared local rule
$\phi$'' do \emph{not} imply equivariance and do \emph{not} support the
irreversibility theorem, because $\eta$ may depend on the cell
arbitrarily. 

Here is a counterexample for state-valued $\eta$ ($d=1$, $n=2$,
$\Sigma=\{0,1\}$, $b=1$): let $\eta_0(\ssb)=(s_0)$, $\eta_1(\ssb)=(0)$, and
$\phi=\mathrm{id}$. All cells share $\phi$, yet the shift-symmetric
configuration $(1,1)$ maps to $(1,0)$, which is not shift-symmetric. 

If one insists that $\eta$ only select \emph{locations} (each cell reads
actual cells, just not translation-uniform ones), equivariance can still
fail on every $\A$, and some location-selecting map breaks irreversibility
exactly when $|\A|=n^d$ is not prime (on $\ZZ_p$, $p$ prime, every
shift-symmetric configuration is constant and is mapped to a constant
one). In $d=1$ the smallest case is $n=4$: with $\phi=\mathrm{id}$, let every cell
read itself except that cell $2$ reads cell $1$; then $(0,1,0,1)$, fixed
by the shift $2$, maps to $(0,1,1,1)$, which is not. 

The original paper by Banda et al.~\cite{banda2019} considers the two-dimensional
case and avoids this pitfall only inside its proofs, where the neighborhood is explicitly
assumed to be given by shared relative offsets \cite[Lemma~A.5]{banda2019}. Our
Definition~\ref{def:ca} in this paper makes the hypothesis explicit from the start.
\end{remark}

\section{Preliminaries}\label{sec:prelim}

In this section, we review some preliminary results about how groups of translations act on the set of configurations.

Throughout, we refer to the prime factorization $n=\prod_{i=1}^{k}p_i^{\alpha_i}$, with $p_1<\dots<p_k$ prime and integers $\alpha_1,\dots,\alpha_k\ge1$. We note that the number of distinct prime factors of $n$ is $k=\omega(n)$. Also recall that, for any prime $p$, the Gaussian binomial $\gaussbinom dr_p=\prod_{i=0}^{r-1}
\frac{p^{d-i}-1}{p^{i+1}-1}$ gives the number of
$r$-dimensional subspaces of $\FF_p^d$, which is the $d$-dimensional vector space over $\FF_p$, the finite field with $p$ elements.

\begin{lemma}[Fixed sets and generated subgroups]\label{lem:fix}
For any given subgroup $H\le\A$ and any given finite family of subgroups $H_1,\dots,H_\ell\le\A$:
\begin{enumerate}[label=(\roman*)]
\item $|\Fix(H)|=\sigma^{|\A|/|H|}$;
\item $\bigcap_{a=1}^{\ell}\Fix(H_a)
      =\Fix\bigl(\langle H_1\cup\dots\cup H_\ell\rangle\bigr)$.
\end{enumerate}
\end{lemma}

\begin{proof}
(i) A configuration $\ssb\in\Fix(H)$ iff $\ssb$ is constant on the cosets of $H$. Translation by $H$ on $\A$ is a free action, so its orbits are the $|\A|/|H|$
cosets, giving a bijection $\Fix(H)\cong\Sigma^{\A/H}$.

(ii) A configuration $\ssb$ lies in the left side iff $H_1\cup\dots\cup
H_\ell\subseteq\Stab(\ssb)$. This happens iff the generated subgroup is contained in
$\Stab(\ssb)$, since $\Stab(\ssb)$ is a subgroup.
\end{proof}

\begin{lemma}[Order of a translation]\label{lem:order}
For $\mathbf 0\neq\vv=(l_1,\dots,l_d)\in\A$, the order of $\vv$ as an element of the additive group $\A$ is given by
\[|\langle\vv\rangle| = \lcm_i\bigl(n/\gcd(l_i,n)\bigr) =
n/\gcd(l_1,\dots,l_d,n),\] and
\[|S_{n^d}(\vv)|=\sigma^{\,n^d/|\langle\vv\rangle|}
=\sigma^{\,n^{d-1}\gcd(l_1,\dots,l_d,n)}.\]
\end{lemma}

\begin{proof}
The order of $\vv$ is the least integer $t\ge1$ with $t \cdot l_i\equiv0\pmod n$ for all
$i$. So $|\langle\vv\rangle|$ must equal $\lcm_i(n/\gcd(l_i,n))$. Comparing $p$-adic valuations,
$v_p(\lcm_i n/\gcd(l_i,n)) = v_p(n)-\min_i\min(v_p(l_i),v_p(n)) =
v_p(n/\gcd(l_1,\dots,l_d,n))$ for every prime $p\mid n$. The count is the
case $H=\langle\vv\rangle$ of Lemma~\ref{lem:fix}(i) above.
\end{proof}

\begin{lemma}[Torsion and lines]\label{lem:torsion}
Let $p\mid n$ be prime. The following hold.

\begin{enumerate}[label=(\roman*)]
\item The $p$-torsion $\A[p]=\{\mathbf x:p\mathbf x=\mathbf 0\}$ is given by
\[\A[p]=\bigl(\tfrac np\ZZ/n\ZZ\bigr)^d.\]

\item The map $\mathbf x=\tfrac np\,\mathbf
y\mapsto\mathbf y\bmod p$ is an isomorphism $\A[p]\cong\FF_p^d$ of
$\FF_p$-vector spaces. 

\item Under this isomorphism, 

$\bullet \;$ elements of order $p$ correspond to
nonzero vectors; 

$\bullet \;$ cyclic subgroups of order $p$ correspond to lines, and
their number is $\gaussbinom d1_p=(p^d-1)/(p-1)$; 

$\bullet \;$ for a set $J$ of
order-$p$ cyclic subgroups, \[|\langle\bigcup J\rangle|=p^{\rk J},\] where
$\rk J$ is the dimension of the span of the corresponding lines.
\end{enumerate}
\end{lemma}

\begin{proof}
Componentwise, $px=0$ in $\ZZ/n\ZZ$ iff $x\in\tfrac np\ZZ/n\ZZ$. The group $\A[p]$
is an elementary abelian $p$-group, and therefore it is an $\FF_p$-vector space. The subgroups of 
$\A[p]$ are its subspaces, and a subgroup generated by a union of lines
is their span.
\end{proof}

\begin{lemma}[Completeness]\label{lem:completeness} For any prime $p$, let $\mathcal L_p$ be the set of cyclic subgroups of order $p$ of $\A$ and let
$\mathcal L=\bigcup_{p\mid n}\mathcal L_p$. Then
\[S_{n^d}=\bigcup_{C\in\mathcal L}\Fix(C).\]
\end{lemma}

\begin{proof}
Since $S_{n^d}$ contains $\Fix(H)$ for \emph{all} nontrivial cyclic subgroups $H$ of $\A$, the direction $\supseteq$ is clear. For the other direction $\subseteq$, we argue as follows. If a configuration $\ssb$ has a nonzero $\ww \in \Stab(\ssb)$, let $o>1$ be its order. Then $o$ divides the
exponent $n$ of $\A$, so any prime divisor $p$ of $o$ must also divide $n$, and
$(o/p)\ww$ has order exactly $p$. It follows that the cyclic subgroup
$C:=\langle(o/p)\ww\rangle\in\mathcal L_p$ satisfies $C\le\Stab(\ssb)$, and so $\ssb\in\Fix(C)$.
\end{proof}

\begin{lemma}[Multi-prime generated subgroups]\label{lem:joins}
Let $J=\bigsqcup_{p\mid n}J_p$ with $J_p\subseteq\mathcal L_p$. Then
\[\bigl|\langle\bigcup J\rangle\bigr|=\prod_{p\mid n}p^{\,\rk J_p}.\]
\end{lemma}

\begin{proof}
Write $\A=\bigoplus_{p\mid n}\A_p$ (primary decomposition). Each
$C\in\mathcal L_p$ lies in $\A[p]\subseteq\A_p$, so $\langle\bigcup
J\rangle=\bigoplus_p\langle\bigcup J_p\rangle$, and each factor has order
$p^{\rk J_p}$ by Lemma~\ref{lem:torsion}.
\end{proof}

\section{Enumeration}\label{sec:enumeration}

In this section, we apply the results above to count the number of shift-symmetric configurations in $S_{n^d}$ for arbitrary finite dimension $d$.

\subsection{The signed spanning-set count}

We begin by introducing a function to denote an alternating sum that will help us to study the various sets of lines that span a finite vector space.

\begin{lemma}\label{lem:mobius}
Let $p$ be prime, $d\ge1$, and $V\le\FF_p^d$ a subspace of dimension $r$,
with $\mathrm{lines}(V)$ its set of $1$-dimensional subspaces. Define 
\[f(V)\;:=\!\!\sum_{\substack{J\subseteq\mathrm{lines}(V)\\
  \operatorname{span}J=V}}\!\!(-1)^{|J|}
\]
Then
\[
  f(V)\;=\;\mu(0,V)\;=\;(-1)^r p^{\binom r2},
\]
where $\mu$ is the M\"obius function of the subspace lattice
\cite[Ex.~3.10.2]{stanley2012}.
\end{lemma}

\begin{proof}
Every $J\subseteq\mathrm{lines}(V)$ spans exactly one $W\le V$, so
\[\sum_{W\le V}f(W)=\sum_{J\subseteq\mathrm{lines}(V)}(-1)^{|J|}
=(1-1)^{|\mathrm{lines}(V)|}=\delta_{V,0},\] since
$\mathrm{lines}(V)\neq\emptyset$ exactly when $V\neq0$. By M\"obius inversion
on the finite subspace lattice, we have $f(V)=\mu(0,V)$. (For more on inversion, see \cite[Prop.~3.7.1]{stanley2012} and \cite[\S4.1.2]{kerber1999}). 

Alternatively, an argument using the $q$-binomial theorem \cite[Eq.~(1.87)]{stanley2012} can be given as follows.  The function mapping a vector space $V$ to
$(-1)^{\dim V}p^{\binom{\dim V}2}$ satisfies the same recursion as $f$ above,
because, by the $q$-binomial theorem  at $x=-1$, we have
\[\sum_{r=0}^{R}\gaussbinom Rr_p(-1)^r p^{\binom r2}
=\prod_{i=0}^{R-1}(1-p^i)=\delta_{R,0}.\]
 This recursion determines $f$ uniquely by induction on $\dim V$. 
\end{proof}

Using the lemma above, we can readily evaluate a similar sum that is taken over all of the sets of lines in $\FF_p^d$ whose span has a given rank $r$.

\begin{lemma}\label{lem:mobius2}
Let $p$ be prime, $d\ge1$, and define 
\[
  c_r(p)\;:=\!\!\sum_{\substack{J\subseteq\mathrm{lines}(\FF_p^d)\\
  \rk J=r}}\!\!(-1)^{|J|}
\]
Then
\begin{equation}\label{eq:cr}
  c_r(p) 
  \;=\;(-1)^r p^{\binom r2}\gaussbinom dr_p ,\qquad 0\le r\le d .
\end{equation}
\end{lemma}

\begin{proof}
Grouping the sets $J$ with $\rk J=r$ by their span, an $r$-dimensional subspace $W$, and
summing $f(W)=(-1)^r p^{\binom r2}$ (Lemma~\ref{lem:mobius}) over the $\gaussbinom dr_p$ subspaces of dimension $r$ gives
\eqref{eq:cr}.
\end{proof}

\subsection{Main theorem}

We are now ready to prove our main result. We begin by defining the coefficient function that appears in our expression for the number of shift-symmetric configurations in $S_{n^d}$ for arbitrary finite dimension $d$.

\begin{definition}[Coefficient function]\label{def:coef}
For a prime $p$ and integers $0\le v\le d$, let
\[r_d(p,v)=p^{\binom v2}\gaussbinom dv_p,\] with $r_d(p,0)=1$. 
\end{definition}

\noindent{\bf Remark.}
    In the case when $d=2$, we observe that $r_2(p,v)=1, p+1, p$ for $v=0,1,2$, which matches the coefficients given in \cite[Thm.~E.19]{banda2019}.

\begin{theorem}[Enumeration in $d$ dimensions]\label{thm:main-enum}
For all $n\ge2$, $d\ge1$, and $\sigma\ge1$:
\begin{equation}\label{eq:main-formula}
    \boxed{\;|S_{n^d}| \;=\; \sum_{\substack{\vv=(v_1,\dots,v_k)\in
    \{0,\dots,d\}^k\\ \vv\neq\mathbf 0}}
    (-1)^{1+|\vv|}\,\prod_{i=1}^{k} r_d(p_i, v_i)\;
    \sigma^{\,n^d\prod_i p_i^{-v_i}}\;}
\end{equation}
where $|\vv|=\sum_i v_i$. All exponents are integers.
\end{theorem}

\begin{proof}
By Lemma~\ref{lem:completeness} and inclusion--exclusion, using
Lemma~\ref{lem:fix} and then Lemma~\ref{lem:joins},
\[
  |S_{n^d}| \;=\; \sum_{\emptyset\neq J\subseteq\mathcal L}(-1)^{|J|+1}
  \bigl|\Fix(\langle\textstyle\bigcup J\rangle)\bigr|
  \;=\; \sum_{\emptyset\neq J}(-1)^{|J|+1}\,
  \sigma^{\,N/\prod_i p_i^{\rk J_{p_i}}} .
\]
The exponent of $\sigma$ depends only on the rank vector
$\vv(J)=(\rk J_{p_1},\dots,\rk J_{p_k})$; the sign is handled by the
factorization below. Adding and subtracting a term $\sigma^N$ for the empty set, we obtain
\[
  |S_{n^d}| \;=\; -\Bigl[\;\sum_{\text{all }J}(-1)^{|J|}\,
  \sigma^{N/\prod p_i^{v_i(J)}} \;-\;\sigma^{N}\Bigr].
\]
Since $(-1)^{|J|}=\prod_i(-1)^{|J_{p_i}|}$ and the parts $J_{p_i}$ vary
independently, the sum factors over rank vectors, with per-prime factors
$c_{v_i}(p_i)$ computed inside $\FF_{p_i}^d$ via Lemma~\ref{lem:torsion}
and Lemma~\ref{lem:mobius2} (note $c_0=1$):
\[
  \sum_{\text{all }J}(-1)^{|J|}\,\sigma^{N/\prod p_i^{v_i}}
  \;=\;\sum_{\vv\in\{0,\dots,d\}^k}\;\prod_{i=1}^k c_{v_i}(p_i)\;
  \sigma^{N/\prod_i p_i^{v_i}} .
\]
The $\vv=\mathbf 0$ term is $\sigma^N$ and cancels, and
$-\prod_i c_{v_i}(p_i)=(-1)^{1+|\vv|}\prod_i r_d(p_i,v_i)$, giving
\eqref{eq:main-formula}. As for the integrality of the exponents, note that $v_i\le d$ and $p_i\mid n$ imply
$\prod_i p_i^{v_i}\mid n^d$.
\end{proof}

\begin{corollary}[Prime lattices]\label{cor:prime}
Fix any prime $n=p$ and integers $d\ge1$ and $\sigma\ge1$. Then
\[|S_{p^d}| = \sum_{r=1}^{d}(-1)^{r+1}\,p^{\binom r2}\gaussbinom dr_p\,
\sigma^{\,p^{d-r}}.\] In particular, 
\[|S_p|=\sigma \quad \mbox{($d=1$),}\]
\[|S_{p^2}|=(p{+}1)\sigma^p-p\sigma \quad \mbox{($d=2$), and}\]
\[|S_{p^3}|=(p^2{+}p{+}1)\sigma^{p^2}-p(p^2{+}p{+}1)\sigma^{p}
+p^3\sigma \quad \mbox{($d=3$).}\]
\end{corollary}

\begin{remark}\label{rem:boundary} We observe the following attributions and boundary checks.

(i) In complement form ($\mathrm{ac}=\sigma^N-|S|$),
Corollary~\ref{cor:prime} is Lemma~7(iv) of Castillo-Ramirez and Gadouleau
\cite{cg2017}; our multi-prime formula \eqref{eq:main-formula} of Theorem~\ref{thm:main-enum} does not
appear there, but we will demonstrate how it can be derived from their framework in
Section~\ref{sec:second-derivation} below.

(ii) In a boundary case where $n=1$, both sides vanish. In a boundary case where $\sigma=1$, the per-prime identity $\sum_{r=0}^d c_r(p)=0$ correctly makes the right side equal to $1$. 

(iii) The case for $d=1$ recovers the classical
count of non-primitive words (proper powers $u^j$, $j\ge2$; equivalently,
words, i.e.\ configurations on $\ZZ_n$, fixed by a nontrivial cyclic
rotation). This is the complement $\sigma^n-\sum_{m\mid n}\mu(m)\,\sigma^{n/m}$
($\mu$ the number-theoretic M\"obius function) of the number of
primitive words; the latter number is $n$ times the Moreau--Witt count of
primitive necklaces \cite{moreau1872,witt1937}.

(iv) The case for $d=2$
recovers \cite[Thm.~E.19]{banda2019} exactly; for prime $n$ it is
\cite[Cor.~E.20]{banda2019}.

(v) For $d=2$, $\sigma=2$, prime $n$, the complement of
Corollary~\ref{cor:prime} appears as an unproved formula line in OEIS
A323863 \cite{oeisA323863}.
\end{remark}

\subsection{Alternative derivation via the subgroup lattice}
\label{sec:second-derivation}

Now we describe how an enumeration matching the above result can also be derived using an alternative line of reasoning.

In \cite[Lemma~7(i)]{cg2017}, Castillo-Ramirez and Gadouleau show that, for any
finite group $G$ of order $N\ge2$ (for $N=1$ the formula is trivial) and $|A|=q\ge2$, the number of aperiodic
configurations is
$\mathrm{ac}(G;A)=\sum_{K\le G}\mu(\{e\},K)\,q^{N/|K|}$, where $\mu$ denotes the
M\"obius function of the subgroup lattice. For $G=\A$, every subgroup
splits along the primary decomposition ($K=\bigoplus_p(K\cap\A_p)$: the
primary components of any $\mathbf x\in K$ are integer multiples of
$\mathbf x$ by CRT, hence lie in $K$), the subgroup lattice is the direct
product of the primary lattices, and $\mu$ is multiplicative across it
\cite[Prop.~3.8.2]{stanley2012}. On each $p$-group factor
$\A_p\cong(\ZZ/p^{\alpha}\ZZ)^d$, by the theorem of Hall
\cite{hall1933} (see \cite[Cor.~3.5]{hio1989}), $\mu(\{e\},K_p)$ vanishes
unless $K_p$ is elementary abelian --- hence contained in the $p$-torsion
$\cong\FF_p^d$ --- in which case
$\mu(\{e\},K_p)=(-1)^{r}p^{\binom r2}$ for rank $r$, and there are
$\gaussbinom dr_p$ such subgroups. Summing
$|S_{n^d}|=\sigma^N-\mathrm{ac}(\A;\Sigma)$ by rank vectors reproduces
\eqref{eq:main-formula}. This second derivation is logically independent
of Lemmas \ref{lem:completeness}--\ref{lem:mobius2}.

\subsection{Computational verification}\label{sec:compver}

Formula \eqref{eq:main-formula} was verified by exact integer computation
against (L1) definition-level brute force over all $\sigma^{n^d}$
configurations and all nonzero shifts, on 29 parameter triples (all with
$n^d\le16$), and (L2) an independent
inclusion--exclusion over all subsets of brute-force-discovered prime-order
subgroups with BFS subgroup closures, on these and 15 further triples
(44 in all), covering
$d\le4$, $n\le60$ (including the multi-prime $n=6,10,12,15,30,60$ and
prime powers $4,8,9$), and $\sigma\le4$; reference values ($\sigma=2$) include
\[
|S_{2^2}|=8,\quad |S_{3^2}|=26,\quad |S_{5^2}|=182,\quad
|S_{2^3}|=72,\quad |S_{3^3}|=6398,\quad |S_{6^2}|=801064 .
\]

\begin{table}[ht]
\centering
\caption{Coefficients $r_d(p,v)=p^{\binom v2}\gaussbinom dv_p$ for small $d$.}\label{tab:coef}
\begin{tabular}{c|c|cccc}
\toprule
$d$ & $v=0$ & $v=1$ & $v=2$ & $v=3$ & $v=4$ \\
\midrule
1 & 1 & $1$ & & & \\
2 & 1 & $p+1$ & $p$ & & \\
3 & 1 & $p^2{+}p{+}1$ & $p(p^2{+}p{+}1)$ & $p^3$ & \\
4 & 1 & $p^3{+}p^2{+}p{+}1$ & $p(p^2{+}1)(p^2{+}p{+}1)$ &
$p^3(p^3{+}p^2{+}p{+}1)$ & $p^6$ \\
\bottomrule
\end{tabular}
\end{table}

\subsection{Finite abelian groups}\label{sec:abelian}

The torus $\ZZ_n^d$ enters the proof of Theorem~\ref{thm:main-enum} only through its order
and its $p$-torsion, which is $\FF_p^d$ for every prime $p\mid n$
(Lemma~\ref{lem:torsion}). For a finite abelian group $G$ of order $N$,
written additively, let $G[p]=\{g\in G: pg=0\}$, an $\FF_p$-vector space, and
let $d_p=\dim_{\FF_p}G[p]$ be the \emph{$p$-rank} of $G$. Configurations
$G\to\Sigma$, translations, stabilizers, and $S_G$ are defined as in
Section~\ref{sec:model} with $\A$ replaced by $G$; a rectangular torus
$\ZZ_{n_1}\times\dots\times\ZZ_{n_m}$ has $N=\prod_i n_i$ and
$d_p=|\{i:p\mid n_i\}|$.

\begin{theorem}[Enumeration on a finite abelian group]\label{thm:abelian}
For every finite abelian group $G$ of order $N$ and every $\sigma\ge1$,
\[
  |S_G| \;=\; \sum_{\substack{\vv=(v_p)_{p\mid N},\ 0\le v_p\le d_p\\
  \vv\neq\mathbf 0}} (-1)^{1+|\vv|}\,\prod_{p\mid N} r_{d_p}(p,v_p)\;
  \sigma^{\,N/\prod_p p^{v_p}} .
\]
For $G=\ZZ_n^d$, where $d_p=d$ for every $p\mid n$, this is
Theorem~\ref{thm:main-enum}.
\end{theorem}
For an abelian $p$-group, Theorem~\ref{thm:abelian} is \cite[Lemma~7(iii)]{cg2017} in complement form,
the elementary abelian subgroups being the subspaces of $G[p]$.

\begin{proof}
The proof of Theorem~\ref{thm:main-enum} applies with $\A$ replaced by $G$
and $\FF_p^d$ by $G[p]\cong\FF_p^{d_p}$: Lemmas~\ref{lem:fix}
and~\ref{lem:completeness} use only that $G$ is a finite abelian group; a
subgroup generated by subgroups of prime order lies in the socle
$\bigoplus_{p\mid N}G[p]$ and is the direct sum of its intersections with
the $G[p]$, which gives Lemma~\ref{lem:joins}; and
Lemmas~\ref{lem:mobius}--\ref{lem:mobius2} apply to each $G[p]$ with $d$
replaced by $d_p$. For a rectangular torus,
$G[p]=\prod_i\ZZ_{n_i}[p]$, and $\ZZ_n[p]\cong\ZZ_p$ if $p\mid n$, while
$\ZZ_n[p]=0$ otherwise.
\end{proof}

By Theorem~\ref{thm:abelian}, $|S_G|$ depends on $G$ only through $N$ and its $p$-ranks; for
$\sigma\ge2$, at fixed $N$ it strictly increases with each $p$-rank, so the $N$ cells arranged as a
one-dimensional ring $\ZZ_N$ carry no more shift-symmetric configurations than any torus with $N$ cells
(Corollary~\ref{cor:abelian-ranks} in Appendix~\ref{app:abelian}, where Remark~\ref{rem:abelian} also extends
Theorem~\ref{thm:kactive} to $G$ and counts the $G$-orbits of aperiodic configurations).

\section{Configurations with $k$ active cells}\label{sec:kactive}

In this section we describe how to refine the enumeration of Theorem~\ref{thm:main-enum} to account for the number of active cells.

To this end, fix any state $a\in\Sigma$ and let $S^a_{n^d,k}$ be the set of
shift-symmetric configurations with exactly $k$ cells in state $a$.
In Theorem~\ref{thm:kactive}, its proof, and
Corollary~\ref{cor:dens} (with $h=h(\vv)$ in Theorem~\ref{thm:kactive},
$h=|H|$ in its proof, and $h=1$ in Corollary~\ref{cor:dens}),
$\binom{N/h}{k/h}:=0$ whenever
$k/h$ is not an integer in $\{0,\dots,N/h\}$, and such a vanishing
binomial annihilates its
whole summand (so the accompanying factor with fractional exponent never
occurs); all other binomial coefficients are the ordinary ones, and
$0^0:=1$ (relevant only at $\sigma=1$, $k=N$). 

In reference \cite{banda2019}, this convention is also implicitly required by Eqs.~(30)--(33)
 (Eq.~(29) needs it only at $\sigma=1$); in
Eqs.~(32)--(33) there, the factors $r(i)$ must be built from the primes
$r_i$ of $\gcd(k,n)$ rather than from the primes $p_i$ of $n$ (as printed,
with the primes in increasing order, both equations give $36$ at $n=6$,
$k=3$, $\sigma=2$, whereas $|S^a_{6^2,3}|=48$).

\begin{theorem}[$k$-active enumeration]\label{thm:kactive}
For all integers $n\ge2$, $d, \sigma\ge1$, and $0\le k\le N=n^d$, with
$h(\vv)=\prod_i p_i^{v_i}$, we have:
\[
  |S^a_{n^d,k}| \;=\; \sum_{\substack{\vv\in\{0,\dots,d\}^{\omega(n)}\\
  \vv\neq\mathbf 0}} (-1)^{1+|\vv|}\,\prod_{i=1}^{\omega(n)} r_d(p_i,v_i)\;
  \binom{\,N/h(\vv)\,}{\,k/h(\vv)\,}\,
  (\sigma-1)^{\,(N-k)/h(\vv)} .
\]
Note: To avoid a notational conflict with the active-cell
count $k$, the rank vectors here run over $\vv\in\{0,\dots,d\}^{\omega(n)}$,
with one coordinate per distinct prime factor.
\end{theorem}

\begin{proof}
For a subgroup $H\le\A$ with $|H|=h$, a configuration in $\Fix(H)$ is a
free labeling of the $N/h$ cosets; it has exactly $k$ cells in state $a$
iff exactly $k/h$ cosets carry the label $a$ --- possible only when
$h$ divides $k$ --- and the remaining $(N-k)/h$ cosets carry labels from
$\Sigma\setminus\{a\}$. Hence
$|\Fix_k(H)|=\binom{N/h}{k/h}(\sigma-1)^{(N-k)/h}$, which is $0$ under the
stated convention unless $h\mid k$. The inclusion--exclusion and
rank-vector grouping of Theorem~\ref{thm:main-enum} apply verbatim, since
the summand again depends on $J$ only through $|\langle\bigcup J\rangle|
= h(\vv(J))$.
\end{proof}

The theorem was verified exactly against definition-level brute force
($n=2$, $d\in\{3,4\}$, $\sigma\le3$, all $k$) and against the independent
subgroup-level evaluator ($n\le6$ at $d=3$; $n\le12$ at $d=2$; all $k$),
including the consistency check $\sum_k|S^a_{n^d,k}|=|S_{n^d}|$.

\section{Bounds}\label{sec:bounds}

We now turn our attention to lower bounds and upper bounds for the number of shift-symmetric configurations, with particular care regarding cases of equality.

\subsection{Lower bound, attained exactly at primes}

For the statement of our lower bound, we will evaluate Gaussian binomials at the (possibly
composite) integer argument $q=n$; since $\gaussbinom dr_q$ is a polynomial
in $q$ with nonnegative integer coefficients, all quantities are integers.

\begin{theorem}[Lower bound; tight iff $n$ prime]\label{thm:lower-bound}
For all integers $n\ge2$, $d\ge1$, $\sigma\ge2$:
\[
    |S_{n^d}| \;\ge\; F_d(n,\sigma)\;:=\;
    \sum_{r=1}^{d}(-1)^{r+1}\,n^{\binom r2}\gaussbinom dr_n\,
    \sigma^{\,n^{d-r}},
\]
with equality if and only if $n$ is prime.
\end{theorem}

The proof of Theorem~\ref{thm:lower-bound}, involving a monotone truncation chain $q^1\le q^2\le\dots\le
q^k=|S_{n^d}|$ over the prime factors, is given in
Appendix~\ref{app:lower}. The argument makes use of an interpretation of the
truncated sums as counts of unions of fixed-point sets on
\emph{rectangular} tori, plus a strict base-case inequality.

For $d=2$, the statement exactly matches
\cite[Lemma~E.21]{banda2019}, and our chain lemma re-proves the
corresponding step of their argument uniformly in $n$ and
$\sigma$.\footnote{The published inequality chains behind the 2D bounds
have small-$n$ gaps (at $n\in\{6,10,14\}$); a companion erratum note
repairs them uniformly. The $d$-dimensional chain of
Appendix~\ref{app:lower} has no exceptional lattice sizes.} A
machine-checked proof of the theorem, taking a different and shorter
route that bypasses the truncation chain altogether, is described in
Section~\ref{sec:formal}.

\begin{remark}[The naive generalization is false]\label{rem:oldbound-false}
A direct generalization of the two-term 2D bound, which would naively predict that
$|S_{n^d}|\ge\sigma^{n^{d-1}}\tfrac{n^d-1}{n-1}
-\sigma\bigl(\tfrac{n^d-1}{n-1}-1\bigr)$ with equality iff $n$ prime, is
\emph{false} for every $d\ge3$. Indeed, at $(n,d,\sigma)=(2,3,2)$ it claims
$|S_{2^3}|\ge100$, while $|S_{2^3}|=72$. For prime $n$, the true count is
the full $d$-term alternating sum of Corollary~\ref{cor:prime}, which
agrees with the naive two-term expression only for $d\le2$ (when $\sigma\ge2$).
\end{remark}

Theorem~\ref{thm:lower-bound} makes precise, and proves, the sense in
which \emph{primes minimize shift-symmetry in every dimension}: the count
$|S_{n^d}|$ meets the universal lower envelope $F_d(n,\sigma)$ exactly at
the primes and exceeds it strictly at every composite $n$.

\subsection{Upper bound}

We now derive two upper bounds on the number of shift-symmetric configurations.

\begin{theorem}[Upper bounds]\label{thm:upper-bound}
For all integers $n\ge2$, $d\ge1$, $\sigma\ge2$:
\begin{enumerate}[label=(\roman*)]
\item (union form)\quad
$\displaystyle |S_{n^d}| \;\le\; \sum_{p\mid n}\gaussbinom d1_p\,
\sigma^{\,n^d/p}$;
\item (closed form)\quad
$\displaystyle |S_{n^d}| \;\le\; \bigl(2^d-1+2\log_2 n\bigr)\,
\sigma^{\,n^d/2}$.
\end{enumerate}
\end{theorem}

\begin{proof}
(i). This is the union bound over Lemma~\ref{lem:completeness}, using
$|\Fix(C)|=\sigma^{N/p}$ for $C\in\mathcal L_p$, as in the proof of
\cite[Thm.~5(ii)]{cg2017}.

\noindent
(ii). If $2\mid n$, the $p=2$ term of (i) is
$(2^d-1)\sigma^{N/2}$ (if $n$ is odd, there is no such term and only
the bound below is needed). We claim every \emph{odd} prime $p\mid n$
contributes at most $2\sigma^{N/2}$; since $n$ has at most $\log_2 n$
prime factors, (ii) follows. The claim, term by term: we must show
$\gaussbinom d1_p\le 2\,\sigma^{\,N(\frac12-\frac1p)}$, and since
$\sigma\ge2$ and $\tfrac12-\tfrac1p\ge\tfrac16$ it suffices that
$\gaussbinom d1_p\le 2^{1+N/6}$.
For $d=1$: $\gaussbinom 11_p=1$. For $d=2$: any odd prime $p\mid n$
satisfies $p\le n$, so $\gaussbinom 21_p=p+1\le n+1$, and
$n+1\le2^{1+n^2/6}$ for every $n\ge3$: directly for $n\in\{3,4,5\}$
($4\le2^{2.5}$, $5\le2^{3.67}$, $6\le2^{5.17}$), and for $n\ge6$ from
$n^2/6\ge n$ and $2^{\,n}\ge n+1$. For $d\ge3$:
$\gaussbinom d1_p\le 2p^{d-1}\le 2n^{d-1}$ and
$N/6\ge n^{d-1}\cdot n/6\ge n^{d-1}/2$, so it suffices to show that
$2^{\,n^{d-1}/2}\ge n^{d-1}$, which is true because $x=n^{d-1}\ge9$ and
$2^{x/2}\ge x$ for $x\ge4$.
\end{proof}

Both parts were verified exactly for $n\le60$ and $d, \sigma\le4$.
For comparison, on $G=\ZZ_n^d$ the general-group bounds give
$|S|\le\sigma^{N-1}$ \cite[Cor.~1.7.2]{gjs2016} and
$|S|\le(N-1)\sigma^{N/p_1}$ \cite[Thm.~5(ii)]{cg2017}, $p_1$ the smallest
prime factor. Bound (i) implies
$|S|\le2\log_2(n)\,n^{d-1}\,\sigma^{N/p_1}$
(Proposition~\ref{prop:rate}), sharpening the prefactor $N-1=n^d-1$ of
\cite{cg2017} at the same exponent for every $n\ge7$ (and for $n\ge5$ when
$d\ge2$; at $n=2$ the union prefactor
$2^d-1$ exactly
ties $N-1$, and for $n\in\{3,4\}$ the two are comparable); for even
$n\ge4$ and $d\ge2$, the closed form (ii) improves the prefactor $2\log_2(n)\,n^{d-1}$
further to $2^d-1+2\log_2 n$. (For odd $n$ the exponent $N/p_1\le N/3$ of (i) is
stronger than the $N/2$ of (ii); the two forms serve different purposes.)
For $d=2$ and $n\ge4$, (ii) also slightly sharpens the
$6\log_2(n)\sigma^{n^2/2}$ of \cite[Lemma~E.22]{banda2019}.

\section{Probability of shift-symmetry}\label{sec:probability}

In this section, we turn from questions of enumeration to the closely-related questions of probability.
How likely is it for randomly selected configurations to possess non-trivial shift-symmetry?

\begin{corollary}[Uniform distribution]\label{cor:prob}
For a configuration drawn uniformly from $\Sigma^{n^d}$, the probability
$P^{\mathrm{unif}}_{n^d}=|S_{n^d}|/\sigma^{n^d}$ satisfies
\[
    F_d(n,\sigma)\,\sigma^{-n^d}
    \;\le\; P^{\mathrm{unif}}_{n^d} \;\le\;
    \bigl(2^d-1+2\log_2 n\bigr)\,\sigma^{-n^d/2},
\]
where $F_d(n, \sigma)$ denotes the lower bound on $|S_{n^d}|$ given in Theorem~\ref{thm:lower-bound}.
So $P^{\mathrm{unif}}_{n^d}$ decays as $\sigma^{-cN}$ with
$N=n^d$ and decay rate
$c=\log_\sigma(1/P^{\mathrm{unif}}_{n^d})/N\in
[\tfrac12-o(1),\,1-\tfrac1{p_1}]$, where $p_1$ is the smallest prime
factor of $n$. The lower end is given by Theorem~\ref{thm:upper-bound}(ii), and the
upper end is the single-subgroup inclusion $\Fix(C)\subseteq S_{n^d}$ for
one order-$p_1$ cyclic subgroup $C$, which gives
$P^{\mathrm{unif}}_{n^d}\ge\sigma^{N/p_1-N}$ (this is also the lower bound
of Proposition~\ref{prop:rate}). Here $o(1)$ refers to any regime with
$N=n^d\to\infty$.
\end{corollary}

\begin{corollary}[Density-uniform distribution]\label{cor:dens}
Let $\sigma\ge2$. Under the density-uniform distribution (first draw
$k\in\{0,\dots,N\}$ uniformly, then a configuration with exactly $k$
cells in state $a$ uniformly), the probability that the drawn configuration lies in $S_{n^d}$ is
\[
    P^{\mathrm{dens}}_{n^d} \;=\; \frac{1}{N+1}\sum_{k=0}^{N}
    \frac{|S^a_{n^d,k}|}{\binom Nk(\sigma-1)^{N-k}},
\]
with $|S^a_{n^d,k}|$ given by Theorem~\ref{thm:kactive}.
\end{corollary}

\begin{figure}[ht]
\centering
\includegraphics[width=\textwidth]{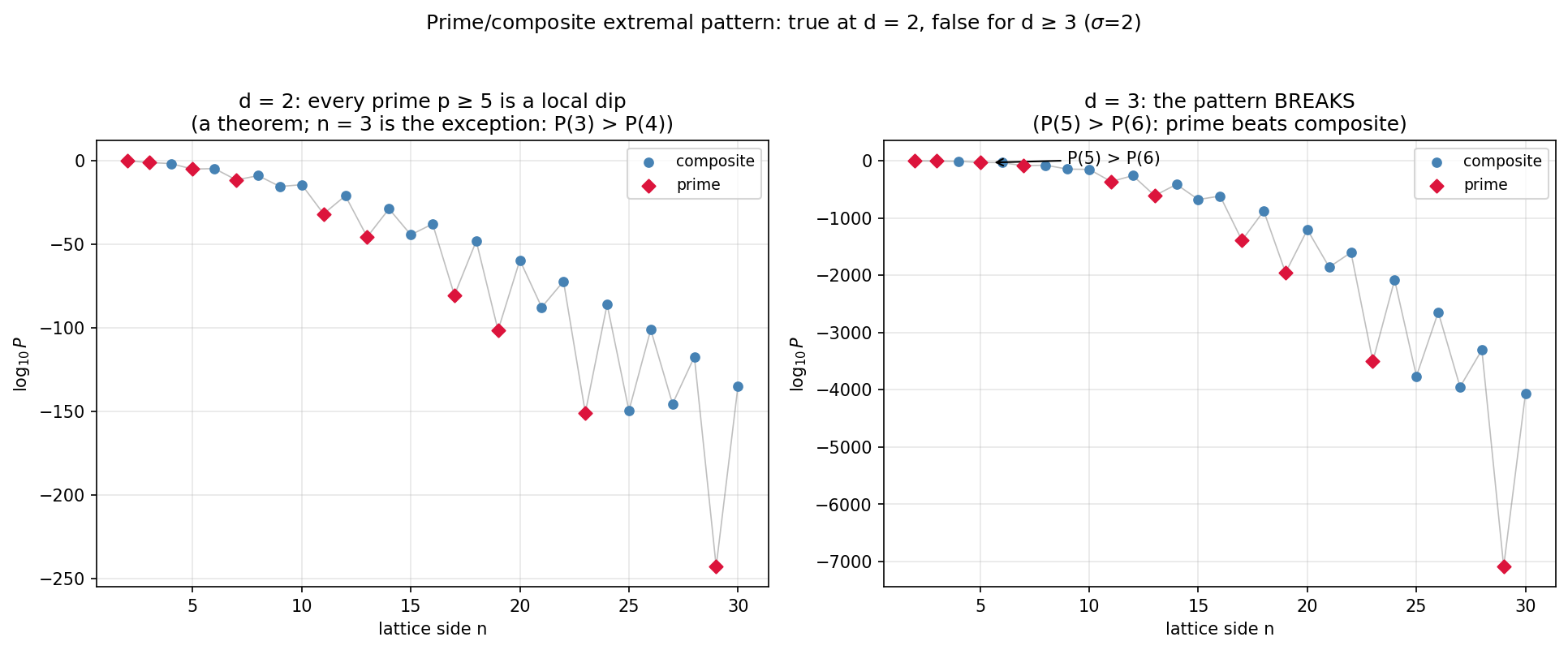}
\caption{The prime/composite extremal pattern of
Observation~\ref{obs:extrema} ($\sigma=2$; the exact formula, evaluated numerically). Left, $d=2$:
$\log_{10}P^{\mathrm{unif}}$ has local dips at every prime
$5\le n\le40$ (shown here for $n\le30$); the boundary case $n=3$ is the
audited exception. Right, $d=3$: the pattern breaks --- already
$P(5^3)>P(6^3)$, so the prime $5$ is no longer a local minimum.}
\label{fig:pattern-flip}
\end{figure}

\begin{observation}[Extremal pattern in $n$: true for $d=2$, false for
$d\ge3$]\label{obs:extrema}
For $d=2$, $\sigma=2$, exact computation for $n\le41$ confirms the pattern
of \cite[Figs.~6,~8]{banda2019}: $P^{\mathrm{unif}}$ has local minima at
every prime $5\le n\le40$ and local maxima at every even $4<n\le40$ (the
boundary is
a genuine exception: $P(3)=26/2^9>736/2^{16}=P(4)$, so $3$ is not a local
minimum). For $d\ge3$ the
pattern \emph{breaks}: at $d=3$, $\sigma=2$,
\[
\log_{10}P(5)=-28.61 \;>\; \log_{10}P(6)=-31.67,
\]
so the prime $5$ is not a local minimum and the composite $6$ is not a
local maximum; likewise the even $n=10$ is not a local maximum
($P^{\mathrm{unif}}_{9^3}>P^{\mathrm{unif}}_{10^3}$), and all three
failures persist at $d=4$.
See Figure~\ref{fig:pattern-flip}.
This is consistent with the inversion law
(Proposition~\ref{prop:inversion}): for the pair $(5,6)$ of the pattern,
$d_c(5)\approx2.58$, and indeed
$P^{\mathrm{unif}}_{5^d}>P^{\mathrm{unif}}_{6^d}$ for every $d\ge3$
($\sigma=2$), which already breaks the pattern (as does
$P^{\mathrm{unif}}_{9^3}>P^{\mathrm{unif}}_{10^3}$), while at $d=3$
($\sigma=2$) no pair $(p,p\pm1)$ with $7\le p\le23$ is inverted. Apart
from two statements that hold for all $d\ge3$ ($\sigma=2$), resp.\ for
all primes $p\ge5$ and all $\sigma\ge2$ --- the inequality
$P^{\mathrm{unif}}_{5^d}>P^{\mathrm{unif}}_{6^d}$ above, and, at $d=2$,
$P^{\mathrm{unif}}_{p^2}<P^{\mathrm{unif}}_{(p\pm1)^2}$
(Observation~\ref{obs:inversion}) --- these are exact finite computations
(all machine-checked, Section~\ref{sec:formal}), not asymptotic claims;
we know of no proof of the even-peak pattern as a theorem for general $n$.
\end{observation}

\section{Dynamics: irreversibility, insolvability, rank monotonicity}
\label{sec:dynamics}

In this section we establish the dynamical consequences of shift-symmetry in configurations.

\begin{theorem}[Irreversibility]\label{thm:irrev}
Let $\Phi$ be translation-equivariant (Definition~\ref{def:ca}). Then, for every configuration $\ssb$,
$\Stab(\Phi(\ssb))\supseteq\Stab(\ssb)$;
in particular, $\ssb\in S_{n^d}(\vv)$ implies
$\Phi(\ssb)\in S_{n^d}(\vv)$, and $\Phi(S_{n^d})\subseteq S_{n^d}$.
\end{theorem}

\begin{proof}
If $\tau_\vv\ssb=\ssb$, then
$\tau_\vv\Phi(\ssb)=\Phi(\tau_\vv\ssb)=\Phi(\ssb)$.
\end{proof}

For cellular automata on any finite group, this inclusion appears in the
proof of Lemma~10 of \cite{cg2017}, whose statement gives the converse
for distinct orbits; see also \cite[Thm.~2(iii)]{cgvn2017}. By Remark~\ref{rem:equivariance-needed},
the equivariance hypothesis cannot be weakened to the assumption ``all cells
share the local rule.''

\begin{corollary}[Insolvability]\label{cor:insolv}
Let $\sigma\ge2$ and $N>1$. Any task requiring, on some input, an output
configuration outside $S_{n^d}$ reachable from an input inside $S_{n^d}$
is unsolvable by translation-equivariant CA on $\ZZ_n^d$. Indeed, symmetric inputs
exist (all constants lie in $S_{n^d}$), and by Theorem~\ref{thm:irrev} no
trajectory leaves $S_{n^d}$. This covers leader election
(the target configuration with a unique distinguished cell is asymmetric
for $N>1$); $d$-dimensional
leader-election algorithms as in \cite{nichitiu1999,stratmann2002}
therefore necessarily
break one of the hypotheses (uniformity, determinism, synchrony) or
restrict the admissible inputs, e.g.\ to asymmetric configurations.
More generally, since $\Stab$ can only grow along a trajectory
(Theorem~\ref{thm:irrev}), every pattern-recognition or image-processing
task whose required output breaks a symmetry of some admissible input is
unsolvable, as in the two-dimensional discussion of \cite{banda2019}.
\end{corollary}

\begin{proposition}[Per-prime rank monotonicity]\label{prop:rank}
For a prime $p\mid n$ define the \emph{$p$-rank} of $\ssb$ as
$\rk_p(\ssb)=\dim_{\FF_p}\bigl(\Stab(\ssb)\cap\A[p]\bigr)\in\{0,\dots,d\}$
(via Lemma~\ref{lem:torsion}). Under any translation-equivariant $\Phi$,
every $p$-rank is non-decreasing along every trajectory:
$\rk_p(\Phi(\ssb))\ge\rk_p(\ssb)$ for all $p\mid n$ simultaneously.
\end{proposition}

\begin{proof}
By Theorem~\ref{thm:irrev}, $\Stab(\Phi(\ssb))\cap\A[p]\supseteq
\Stab(\ssb)\cap\A[p]$, and both are $\FF_p$-subspaces.
\end{proof}

For composite $n$ there is one rank per prime factor, and all of them are
monotone. Geometrically, $\rk_p=1$ configurations are
periodic along one direction of the $p$-torsion (``stripes''/``slabs''),
$\rk_p=d$ configurations are exactly those fixed by the whole $p$-torsion
$\A[p]=(\tfrac np\ZZ/n\ZZ)^d$; each superlevel set $\{\rk_p\ge r\}$ is
absorbing (no trajectory leaves it) by Proposition~\ref{prop:rank}. Kornyak's
observation that deterministic dynamics cannot decrease symmetry
\cite{kornyak2010} is the general heuristic; the proposition is its exact
translation-stabilizer form.

\section{Detection}\label{sec:detection}

In this section, we discuss how to determine whether configurations possess shift-symmetry.

The detector given here generalizes that given in an earlier paper \cite[Sec.~V]{banda2019}. Assume the prime
factorization of $n$ is given (otherwise add one factorization of an
integer of $\log_2 n$ bits), and count cell-read operations; $d$ is part
of the input.

\begin{quote}
\textbf{Algorithm.} For each prime $p\mid n$ and each of the
$\gaussbinom d1_p$ order-$p$ cyclic subgroups $C=\langle\ww\rangle$
(enumerated via one representative line each,
Lemma~\ref{lem:torsion}): traverse the cells in a fixed order (which may depend on the subgroup); for each
yet-unvisited cell, walk its $C$-orbit ($p$ cells) checking equality;
abort the subgroup at the first failing orbit. Report ``symmetric'' at
the first subgroup all of whose orbits pass; report ``asymmetric'' if all
subgroups fail.
\end{quote}

Correctness is validated by Lemma~\ref{lem:completeness}: a configuration is
shift-symmetric iff it is fixed by some order-$p$ cyclic subgroup with
$p\mid n$.

Throughout this section we charge one read per visited cell: the first
cell of an orbit walk (the \emph{anchor}) is read once and its value
cached, and each subsequent cell of the orbit is read once and compared
to the cached value. A failing orbit walk may stop at its first mismatch or be completed; all
results of this section hold in both variants, except the exact factor in
Remark~\ref{rem:polylog}, which refers to completed walks.

In what follows, let $\mathrm{sopf}_d(n)=\sum_{p\mid n}\gaussbinom d1_p
=\sum_{p\mid n}\frac{p^d-1}{p-1}$ denote the total number of order-$p$
cyclic subgroups over all primes $p\mid n$ (the $d$-dimensional analogue
of the integer logarithm).

\begin{theorem}[Worst case]\label{thm:detect-worst}
The algorithm reads at most
$N\cdot\mathrm{sopf}_d(n)=O(\mathrm{sopf}_d(n)\cdot n^d)$ cells.
\end{theorem}

\begin{proof}
Per subgroup, an orbit walk of length $p$ performs $p$ reads, and the
walked orbits are disjoint, so one subgroup costs at most $N$ reads;
there are $\mathrm{sopf}_d(n)$ subgroups.
\end{proof}

\begin{proposition}[Average case, uniform random input]
\label{prop:detect-avg}
For $\ssb$ uniform on $\Sigma^{n^d}$, the test of any single order-$p$
subgroup, run on its own, visits in expectation exactly
$\frac{1-Q^{m}}{1-Q}$ orbits, with $Q=\sigma^{1-p}$ and $m=N/p$. Consequently, the
expected total number of cells read by the full algorithm is at most
\[
    \sum_{p\mid n}\gaussbinom d1_p\;\frac{p}{1-\sigma^{1-p}}
    \;\le\; 2\sum_{p\mid n}p\,\gaussbinom d1_p
    \;=\;O\Bigl(\sum_{p\mid n}p^{\,d}\Bigr).
\]
\end{proposition}

\begin{proof}
Orbits of a fixed subgroup partition the cells, so under the product
measure the events ``orbit is monochromatic'' are independent, each of
probability $\sigma^{1-p}$; the walk visits $i$ orbits iff the first
$i-1$ pass, giving the truncated geometric expectation
$\sum_{i=0}^{m-1}Q^i$. Summing $p\cdot(1-Q)^{-1}$ over all subgroups
bounds the total; $(1-\sigma^{1-p})^{-1}\le2$ for $\sigma,p\ge2$.
\end{proof}

The per-generator expectation was confirmed to match full-enumeration
averages \emph{exactly} for $(n,d,\sigma)\in\{(2,3,2),(2,3,3),(2,4,2)\}$
and, in two dimensions, for $n\le4$ with $\sigma=2$ and for $(n,\sigma)=(3,3)$.

\begin{remark}[On ``polylogarithmic'' average cases]\label{rem:polylog}
$\sum_{p\mid n}p^d$ is \emph{not} polylogarithmic in $N$ in general: for
prime $n$ it equals $n^d=N$. The polylogarithmic figures quoted in
\cite[Sec.~V--VI]{banda2019} for $d=2$ are stated there for
\emph{randomly chosen} $n$ and rest on $\sum_{p\mid n}p$ and
$\sum_{p\mid n}p^2$ being polylogarithmic; but for each fixed dimension $d$,
$\sum_{p\mid n}p^d$ is polylogarithmic in $N$ (with a degree depending on
$d$) when all prime factors of $n$ are polylogarithmically small, whereas for any fixed $C$ and $j$, the $n$
with $\sum_{p\mid n}p^d\le C\log^j N$ have natural density zero, so this
bound is not polylogarithmic for almost all $n$. Note also that, for prime $n$, the average case
$O(N)$ still improves the worst case, which is of exact order
$N\cdot n^{d-1}$, by a factor of at least order $n^{d-1}$ (of exact order
$n^{d-1}$ when every failing orbit walk is completed, as in the proofs of
Theorem~\ref{thm:detect-worst} and Proposition~\ref{prop:detect-avg}).
\end{remark}

\begin{remark}[Read complexity]\label{rem:reads}
With memory for all $N$ cells, $N$ reads suffice: read every cell once and
test in memory. For $N\ge2$ and $\sigma\ge2$, $N$ reads are also necessary
in the worst case: every correct deterministic algorithm, adaptive or not, makes at least $N$
reads on every constant configuration. Indeed, a constant
configuration is shift-symmetric, while the configuration that differs from
it in a single unread cell has trivial stabilizer and yields the same reads,
hence the same answer. Writing $I_a$ for the indicator of the cells in
state $a$, a translation $\ww$ lies in $\Stab(\ssb)$ if and only if
$\sum_{a\in\Sigma}\sum_{\uu\in\A}I_a(\uu)\,I_a(\uu+\ww)=N$.
\end{remark}

If the configuration is held in memory, the sums of Remark~\ref{rem:reads} are obtained for all $\ww$ at once by $d$-dimensional
fast Fourier transforms \cite{cooleytukey1965,bluestein1970}, one per symbol
that occurs in $\ssb$, in $O(\min(\sigma,N)\cdot N\log N)$ arithmetic
operations, without rounding error if a number-theoretic transform is used
\cite{pollard1971}. This improves on the worst case of
Theorem~\ref{thm:detect-worst} when $\mathrm{sopf}_d(n)$ is large compared
with $\min(\sigma,N)\log N$; for fixed $\sigma$, this is the case for
$d\ge2$ as the prime $n$ grows, and for $n=2$ as $d$ grows. For fixed $d$
and $n$ a power of $2$, for $d=1$, and on average
(Proposition~\ref{prop:detect-avg}), the orbit walks are cheaper.

Related decision problems on toroidal arrays --- shift-equivalence of
tilings and recognition of toroidal codes, including undecidability
results --- were recently analyzed by Anselmo, Cavallaro, Madonia, and
Selmi \cite{anselmo2025}; detecting a nontrivial translation stabilizer,
as here, is the tractable member of that family.

\section{Crossover laws and empirical observations in higher dimensions}
\label{sec:empirical}

We first prove two crossover laws that organize all the observations made in this section; the observations after them are computational findings with
exact arithmetic, labeled as such. Where an item after the two laws
states an exact claim for all parameters (such as the composite
advantage at $d=2$ and the deepening of the inversion in
Observation~\ref{obs:inversion}), that claim is a theorem,
machine-checked in Section~\ref{sec:formal}.

\begin{proposition}[Decay-rate law: only the smallest prime factor
survives]\label{prop:rate}
Let $p_1=p_1(n)$ be the smallest prime factor of $n$ and $N=n^d$. For all
$n\ge2$, $d\ge1$, $\sigma\ge2$:
\[
    \sigma^{N/p_1}\;\le\;|S_{n^d}|\;\le\;
    2\log_2(n)\,n^{d-1}\,\sigma^{N/p_1}.
\]
Consequently the decay rate $c(n,d):=-\log_\sigma
P^{\mathrm{unif}}_{n^d}/N$ satisfies
\[
    0\;\le\;\Bigl(1-\frac1{p_1}\Bigr)-c(n,d)\;\le\;
    \frac{\log_\sigma\bigl(2\,n^{d-1}\log_2 n\bigr)}{N}
    \;\xrightarrow[d\to\infty]{}\;0 ,
\]
so $c(n,d)\to1-1/p_1(n)$. Asymptotically in $d$, lattice sides fall into
\emph{universality classes} determined solely by their smallest prime
factor, with even $n$ on the extremal class $c=1/2$ (the floor permitted
by Theorem~\ref{thm:upper-bound}) and primes $p$ on $c=1-1/p$.
\end{proposition}

\begin{proof}
Lower bound: $\Fix(C)\subseteq S_{n^d}$ for any single order-$p_1$ cyclic
subgroup $C$, and $|\Fix(C)|=\sigma^{N/p_1}$ (Lemma~\ref{lem:fix}).
Upper bound: in Theorem~\ref{thm:upper-bound}(i), every exponent satisfies
$N/p\le N/p_1$, and
$\sum_{p\mid n}\gaussbinom d1_p\le\sum_{p\mid n}2p^{d-1}
\le2\,\omega(n)\,n^{d-1}\le2\log_2(n)\,n^{d-1}$.
The prefactor $N-1$ of \cite[Thm.~5(ii)]{cg2017} gives the same limit.
\end{proof}

\begin{figure}[ht]
\centering
\includegraphics[width=\textwidth]{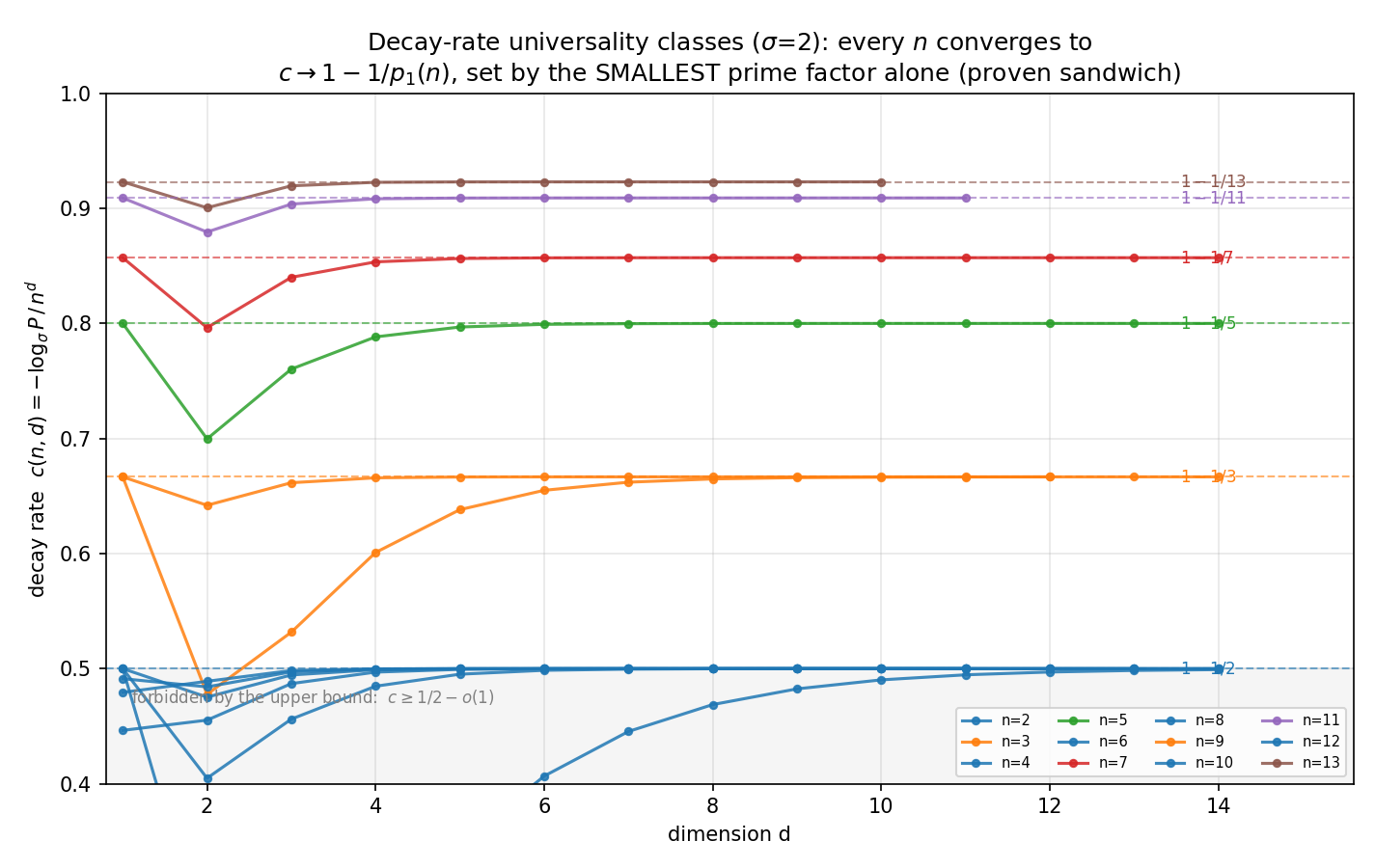}
\caption{Decay-rate universality classes
(Proposition~\ref{prop:rate}): the measured rate
$c(n,d)=-\log_\sigma P^{\mathrm{unif}}_{n^d}/n^d$ at $\sigma=2$ for
$2\le n\le13$, plotted against $d$ (the exact formula, evaluated numerically). Every lattice
side converges to the level $1-1/p_1(n)$ set by its smallest prime
factor alone; all even $n$ collapse onto the extremal class $c=1/2$,
the floor permitted by Theorem~\ref{thm:upper-bound} (shaded region
excluded asymptotically).}
\label{fig:rate-classes}
\end{figure}

Figure~\ref{fig:rate-classes} shows the convergence in numerical evaluations
of the exact formula.

\begin{proposition}[Inversion law]\label{prop:inversion}
Let $p\ge3$ be prime, $\sigma\ge2$, and set
$\varphi(d)=(1-\tfrac1p)\,p^d-\tfrac12(p+1)^d$. Then:
\begin{enumerate}[label=(\roman*)]
\item if $\varphi(d)>\log_\sigma\bigl(2\log_2(p)\,p^{d-1}\bigr)$, then
$P^{\mathrm{unif}}_{(p+1)^d}>P^{\mathrm{unif}}_{p^d}$ (composite
advantage);
\item if $\varphi(d)<-\log_\sigma\bigl(2\log_2(p{+}1)\,(p{+}1)^{d-1}\bigr)$,
then $P^{\mathrm{unif}}_{(p+1)^d}<P^{\mathrm{unif}}_{p^d}$ (prime
advantage).
\end{enumerate}
The function $\varphi$, regarded as a function of real $d$, is positive
below and negative above its unique zero
\[
    d_c(p)\;=\;\frac{\ln\bigl(2-\tfrac2p\bigr)}{\ln\bigl(1+\tfrac1p\bigr)}
    \;=\;p\ln2+O(1)\qquad(p\to\infty),
\]
so the two-dimensional ``composite advantage'' inverts within an
explicitly bounded window around $d_c(p)$, a dimension that grows
\emph{linearly} in the prime: (i) applies to every integer $d\ge1$ with
$d\le d_c(p)-20$, and (ii) to every integer $d\ge d_c(p)+20$.
\end{proposition}

\begin{proof}
(i) By Proposition~\ref{prop:rate} applied to $n=p+1$ (whose smallest
prime factor is $2$) and to $n=p$:
$P_{(p+1)^d}\ge\sigma^{-(p+1)^d/2}$ and
$P_{p^d}\le2\log_2(p)\,p^{d-1}\sigma^{-p^d(1-1/p)}$, so
$\log_\sigma\bigl(P_{(p+1)^d}/P_{p^d}\bigr)\ge
\varphi(d)-\log_\sigma(2\log_2(p)p^{d-1})>0$. (ii) is symmetric.
Monotonicity of $d\mapsto\varphi(d)/p^d=(1-\tfrac1p)-
\tfrac12(1+\tfrac1p)^d$ gives the unique sign change at $d_c(p)$, and
$\ln(2-2/p)/\ln(1+1/p)=p\ln2+O(1)$ by expanding both logarithms. The window bound follows by elementary estimates of $\varphi$ and of the thresholds in (i) and (ii); it is machine-checked (Section~\ref{sec:formal}), and the constant $20$ is not optimized.
\end{proof}

\begin{figure}[ht]
\centering
\includegraphics[width=\textwidth]{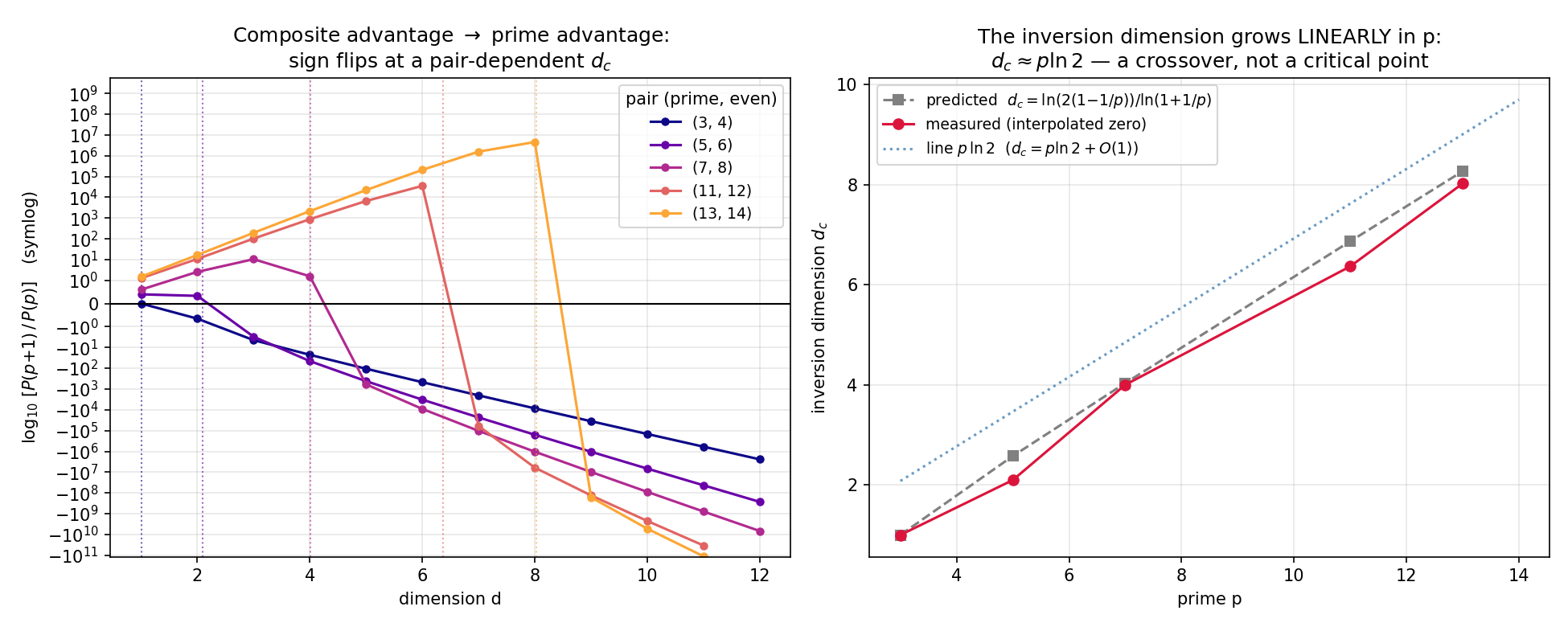}
\caption{The inversion law (Proposition~\ref{prop:inversion},
Remark~\ref{rem:dc-measured}; $\sigma=2$; the exact formula, evaluated numerically). Left: the
ratio $\log_{10}\bigl[P^{\mathrm{unif}}_{(p+1)^d}/P^{\mathrm{unif}}_{p^d}
\bigr]$ for the prime--even pairs $(p,p{+}1)$, $p\in\{3,5,7,11,13\}$,
flips sign at a pair-dependent inversion dimension (symlog scale;
dotted verticals mark the interpolated zeros). Right: the measured
inversion dimensions against the prediction
$d_c(p)=\ln(2-2/p)/\ln(1+1/p)$ and the line $p\ln2$ ($d_c(p)=p\ln2+O(1)$;
the linear asymptote itself is $p\ln2+\tfrac12\ln2-1$).}
\label{fig:inversion}
\end{figure}

\begin{remark}[Measured inversion points]\label{rem:dc-measured}
Linearly interpolated zeros of $d\mapsto\log_{10}[P_{(p+1)^d}/P_{p^d}]$
between the consecutive integer dimensions at its first sign change, in
exact data ($\sigma=2$), against $d_c(p)$:
\begin{center}
\begin{tabular}{lccccc}
\toprule
$p$ & 3 & 5 & 7 & 11 & 13\\
\midrule
measured & $1.00$ & $2.10$ & $4.00$ & $6.37$ & $8.03$\\
$d_c(p)$ & $1.00$ & $2.58$ & $4.04$ & $6.87$ & $8.27$\\
\bottomrule
\end{tabular}
\end{center}
Measured values sit at or below $d_c$, mainly because of the
interpolation itself: interpolating the bare exponent $\varphi$ (without
coefficient corrections) linearly between $\lfloor d_c\rfloor$ and
$\lfloor d_c\rfloor+1$ already gives a zero below $d_c(p)$ for
$p=5,7,11,13$ (at $p=3$ both equal $1$); the coefficient corrections (the factors $|S_{n^d}|/\sigma^{N/p_1}$, which
Proposition~\ref{prop:rate} bounds on both sides) favor the prime side but
account for at most a quarter of the gap (at most $10^{-3}$ of it for
$p=11,13$). At
$p=3$ the boundary case is exact and alphabet-independent:
$P^{\mathrm{unif}}_{3}=P^{\mathrm{unif}}_{4}=\sigma^{-2}$ at $d=1$, since
$|S_3|=\sigma$ and $|S_4|=\sigma^2$. Figure~\ref{fig:inversion} displays
both the sign flips and the measured-versus-predicted inversion
dimensions.
\end{remark}

The remaining observations are computational, except for exact claims
stated for all parameters (such as the composite advantage at $d=2$ and
the deepening of the inversion in Observation~\ref{obs:inversion}),
which are theorems (machine-checked, Section~\ref{sec:formal}).

\begin{observation}[Single-term crossover dimension]\label{obs:dstar}
For prime $n=p$, define $d^*(p;\theta)$ as the least $d\ge2$ at which the
second term of Corollary~\ref{cor:prime} falls below a fraction $\theta$
of the first, i.e.\ $|T_2/T_1|<\theta$ where
$|T_2/T_1| = p\gaussbinom d2_p\gaussbinom d1_p^{-1}
\sigma^{-p^{d-2}(p-1)}$. Computing $d^*$ for all 95 primes up to 499
($\sigma=2$): only $p\in\{2,3,5\}$ have $d^*>2$ at $\theta=0.05$; a fit
$d^*\approx1+C/\ln p$ to a fractional version of $d^*$ (interpolating
$\log|T_2/T_1|$ linearly between consecutive integer dimensions) gives
$C\approx2.2$--$2.65$ depending on $\theta$
(about $2.5\pm0.3$), but with only two to four unclamped data points per
threshold (for the integer $d^*$: four, three, and two at
$\theta=0.01,0.05,0.10$) and
systematic residuals this is an indicative, threshold-dependent summary,
not an estimate of a universal constant. The driver is the
doubly-exponential suppression $\sigma^{-p^{d-2}(p-1)}$.
\end{observation}

\begin{observation}[Prime--composite gap inversion]\label{obs:inversion}
At $d=2$, for every $\sigma\ge2$, every prime $p\ge5$ has smaller
$P^{\mathrm{unif}}$ than both adjacent integers $p\pm1$ (a theorem,
machine-checked in Section~\ref{sec:formal}; the pair $(3,4)$
is the exception: $P^{\mathrm{unif}}_{4^2}<P^{\mathrm{unif}}_{3^2}$;
cf.\ Observation~\ref{obs:extrema}); by $d=3$ ($\sigma=2$) this inverts for the pair $(5,6)$
(Observation~\ref{obs:extrema}), but for no pair $(p,p\pm1)$ with
$7\le p\le23$, and the inversion deepens with $d$: for every prime $p\ge3$ and every
$\sigma\ge2$, once $P^{\mathrm{unif}}_{(p+1)^d}<P^{\mathrm{unif}}_{p^d}$,
the ratio $P^{\mathrm{unif}}_{(p+1)^d}/P^{\mathrm{unif}}_{p^d}$ decreases
strictly with every further dimension (a theorem, machine-checked in
Section~\ref{sec:formal}; e.g., $P(6)/P(5)$ falls from ${\sim}10^{0.33}$ at $d=2$ to
$10^{-419.2}$ at $d=5$, $\sigma=2$).
\end{observation}

\begin{remark}[Duality]
$\gaussbinom dr_p=\gaussbinom d{d-r}_p$ for $0\le r\le d$, so $\FF_p^d$
has equally many subspaces of dimensions $r$ and $d-r$ (equivalently, for
a prime $p\mid n$, $\A[p]$ has equally many subgroups of orders $p^r$ and
$p^{d-r}$). We state this purely as the standard
symmetry of Gaussian binomials; we do not attach further structural
claims to it.
\end{remark}

\paragraph{Reproducibility.}
\begin{sloppypar}
All counts and verification computations use exact integers; the figures
evaluate the exact formula numerically in log space (cross-checked against
exact integers where feasible). The verification suite
(definition-level brute force; independent subgroup-level evaluators for
$|S|$ and $|S^a_k|$; bound and extremal-pattern scans; the
detection-expectation enumeration) and the chart-generation scripts
(\texttt{compute\_charts.py}, \texttt{compute\_criticality.py},
\texttt{compute\_C\_exploration.py}, \texttt{compute\_dimension\_insights.py},
which produces Figures~\ref{fig:pattern-flip}--\ref{fig:inversion}, and
\texttt{compute\_example\_configs.py}, which produces Figure~\ref{fig:examples}
and checks its stabilizers and count) are in the project
repository. Commands, environments, outputs, and SHA-256 checksums are
recorded in \texttt{research/\allowbreak experiments/\allowbreak
README.md}.
\end{sloppypar}

\section{Formal verification}\label{sec:formal}

\begin{figure}[t]
\centering
\begin{tikzpicture}[
    execute at begin picture={\hyphenpenalty=10000\relax},
    box/.style={draw, rounded corners=2pt, align=center,
      font=\footnotesize, text width=2.6cm, inner sep=3pt,
      minimum height=1.5cm},
    chk/.style={draw, dashed, rounded corners=2pt, align=center,
      font=\footnotesize, inner sep=3pt},
    arr/.style={-{Stealth[length=2mm]}, thick},
    lab/.style={font=\scriptsize\itshape, above=1pt}]
  \node[box] (paper) {\textbf{Paper}\\ statements in the \TeX{} source and PDF};
  \node[box, right=1.3cm of paper] (frozen) {\textbf{Frozen \texttt{Prop}s}\\
    370 statements, each tied to its passage; fingerprinted};
  \node[box, right=1.3cm of frozen] (proofs) {\textbf{Proofs}\\
    \texttt{X\_proved} of type \texttt{XStatement}};
  \node[box, right=1.3cm of proofs] (kernel) {\textbf{Lean kernel}\\
    build audit 370/370; three standard axioms};
  \draw[arr] (paper) -- node[lab] {formalize} (frozen);
  \draw[arr] (frozen) -- node[lab] {prove} (proofs);
  \draw[arr] (proofs) -- node[lab] {check} (kernel);
  \coordinate (m) at ($(paper.east)!0.5!(frozen.west)$);
  \node[chk, text width=3.4cm, above=0.6cm of paper] (comp)
    {\textbf{Computation}\\ exhaustive enumeration and an independent
    evaluator; runs logged with SHA-256};
  \node[chk, text width=4.6cm, anchor=north] (fid) at ($(m)+(0,-1.2)$)
    {\textbf{Fidelity checks}\\ independent reviews of each group; hand
    audit of 27 items, 529 links (Section~\ref{sec:conclusion})};
  \draw[arr, dashed] (comp) -- (paper);
  \draw[arr, dashed] (fid.north) -- (m);
\end{tikzpicture}
\caption{The verification pipeline. Statements of the paper are frozen as
Lean \texttt{Prop}s (all but five before their proofs were written); each
proof \texttt{X\_proved} must have exactly the frozen type; the Lean kernel
checks the proofs, and an audit module in the build fails unless every
frozen statement (currently $370$) is proved from the three standard axioms.
Whether a frozen statement says what the paper says is checked separately
(dashed), by independent fidelity reviews and, for the theorems, lemmas,
propositions, and corollaries, by the hand audit described in
Section~\ref{sec:conclusion}; formulas and numbers are, in addition, checked
by exact computation.}
\label{fig:pipeline}
\end{figure}
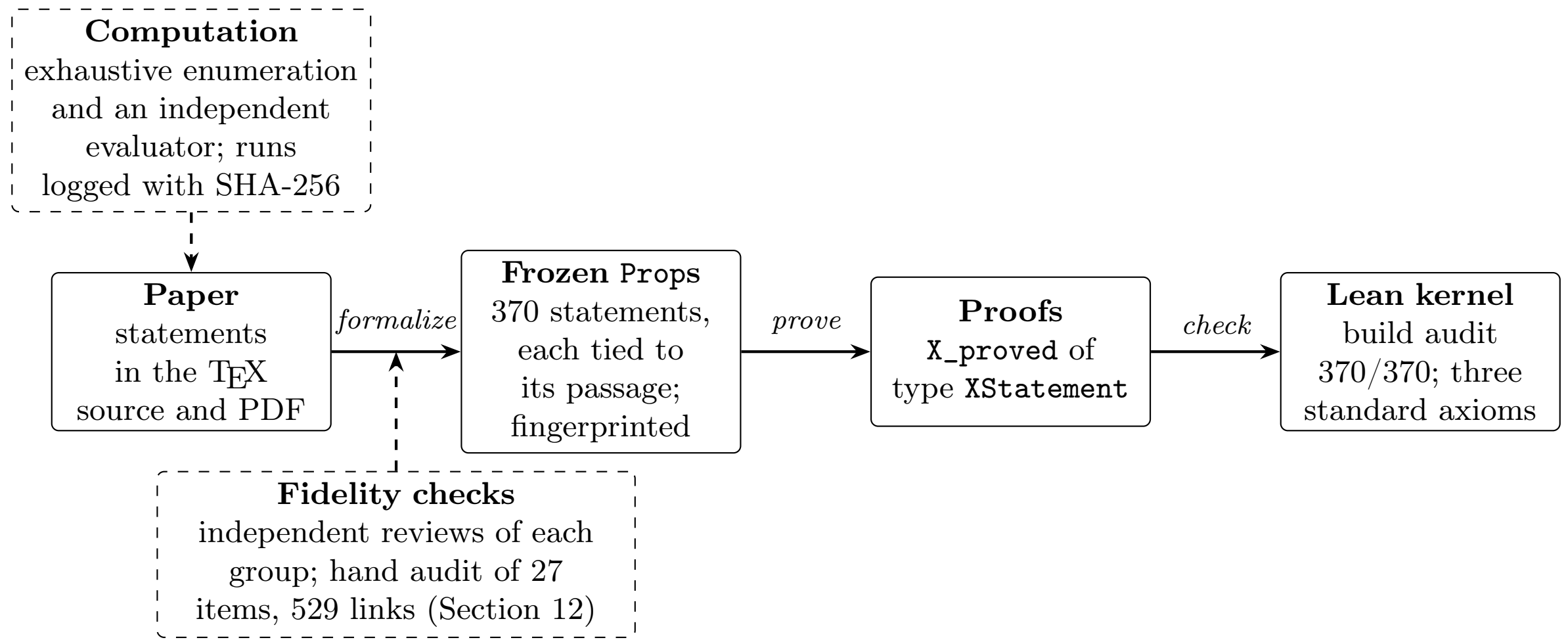

The mathematical content of every numbered statement of
Sections~\ref{sec:model}--\ref{sec:empirical} and
Appendices~\ref{app:lower} and~\ref{app:abelian} has been formally verified in the Lean~4 proof
assistant (v4.33.0-rc1) on top of the \textsf{mathlib} library: that of
all definitions, lemmas, theorems, propositions, corollaries, remarks,
and computational observations there, of Table~1, of the table of
Remark~\ref{rem:dc-measured}, and of the captions of
Figures~\ref{fig:pattern-flip}--\ref{fig:inversion}, together with every
mathematical claim in the surrounding (non-proof) prose that is precise
enough to be stated formally. Statistical summaries (such as the fitted
constant $C$ of Observation~\ref{obs:dstar}), vague, bibliographic, and
record-keeping sentences, comparisons of operation counts that rest on
cited results (after Remark~\ref{rem:reads}), and the summary table of the
concluding section are not formalized. Where the verification found
a mathematical claim ambiguous or false, the text
above states the corrected claim, and that is the claim certified. The
illustrative Figures~\ref{fig:examples} and~\ref{fig:pipeline} and the
applications discussed in Section~\ref{sec:conclusion} are not formalized;
the stabilizers in the caption of Figure~\ref{fig:examples} were checked by
exact computation, and its count is one of the machine-checked reference
values of Section~\ref{sec:compver}.
Figure~\ref{fig:pipeline} summarizes the verification.

The development has two layers. The core layer ($\approx$3200 lines: 15
lemma modules plus definition, statement, experiment, and audit files)
defines the model from first principles --- lattice, configurations,
shifts, stabilizers, fixed sets, equivariant maps --- and proves the two
main theorems:
\begin{itemize}
    \item \texttt{enumeration\_proved}: Theorem~\ref{thm:main-enum} for
    all $n\ge2$, $d\ge1$, $\sigma\ge1$, exactly as stated ($\sigma\ge1$
    being weaker than the standing alphabet assumption $\sigma\ge2$ of
    Section~\ref{sec:model});
    \item \texttt{lowerBound\_proved}: Theorem~\ref{thm:lower-bound}
    verbatim ($n\ge2$, $d\ge1$, $\sigma\ge2$), including both directions
    of the equality-iff-prime characterization.
\end{itemize}
The paper-wide layer ($\approx$27\,000 lines in 67 modules) states and
proves everything else: 370 \texttt{Prop}-valued statements (including
selected proof-internal claims, representation bridges, a
few disclosed strengthenings, and machine-checked evidence for the
corrections), organized by part of the paper ---
Section~\ref{sec:model} with the convention sentence opening
Section~\ref{sec:prelim} (31), the rest of
Section~\ref{sec:prelim} (20), Lemmas~\ref{lem:mobius}--\ref{lem:mobius2},
Definition~\ref{def:coef}, Table~\ref{tab:coef} and the Duality remark (22), the rest of
Section~\ref{sec:enumeration} (48, and 5 more for the necklace
identification in Remark~\ref{rem:boundary}(iii)), Section~\ref{sec:abelian} with Appendix~\ref{app:abelian} (17),
Section~\ref{sec:kactive} with
Corollary~\ref{cor:dens} (31), Section~\ref{sec:bounds} with
Corollary~\ref{cor:prob} (41), Section~\ref{sec:dynamics} with the
caption of Figure~\ref{fig:pattern-flip} (33),
Section~\ref{sec:detection} (31),
Propositions~\ref{prop:rate}--\ref{prop:inversion} with the caption of
Figure~\ref{fig:rate-classes} (24),
Appendix~\ref{app:lower} (33), and
Observation~\ref{obs:extrema}, Remark~\ref{rem:dc-measured} with the
caption of Figure~\ref{fig:inversion}, and
Observations~\ref{obs:dstar}--\ref{obs:inversion} (34). Some of its proofs
also use certified results of a companion development (two further
modules, $\approx$440 lines): the proof of Theorem~\ref{thm:kactive}
uses the relative enumeration theorem, which extends the
inclusion--exclusion behind Theorem~\ref{thm:main-enum} from the set of
all configurations to an arbitrary set of configurations, and proofs in
Appendix~\ref{app:lower} use its weight-generic regrouping lemma.
Multi-part results
are stated part by part and as a whole; asymptotic claims carry explicit
quantifiers; a rounded decimal is stated as the interval its rounding
asserts; the detection algorithm of Section~\ref{sec:detection} is a Lean
function with explicit read counting, and its statements are quantified
over every enumeration of the subgroups, every traversal order, and both
variants of a failing orbit walk (except the exact factor of
Remark~\ref{rem:polylog}, which is stated for completed walks only); the
density-uniform distribution of
Corollary~\ref{cor:dens} is a probability mass function in
\textsf{mathlib}. Numerical claims are proved with exact integer
certificates checked by the Lean kernel, combined, for large lattices
(up to $14^9$ cells in Remark~\ref{rem:dc-measured}), with proved
analytic bounds such as the one-fixed-set and union bounds; no compiled
code or floating-point arithmetic is trusted.
Several of the 370 statements record, as machine-checked facts, why an
earlier wording or a cited formula fails in its literal reading.

All formal statements were frozen as \texttt{Prop}-valued definitions --- all but
five of them \emph{before} their proofs were written --- each with a documentation string
that locates the passage it formalizes, in most cases quoting it
verbatim, and records its representation choices
(e.g.\ the alphabet is realized as $\{0,\dots,\sigma-1\}$, harmless since
only $|\Sigma|$ enters; rank vectors are reindexed by divisors of
$(p_1\cdots p_k)^d$; quotients are exact), and each group of statements
passed independent reviews of its fidelity to the paper. In addition, the
correspondence of every theorem, lemma, proposition, and corollary with its
frozen statement is audited by hand (Section~\ref{sec:conclusion}). Every frozen
statement \texttt{XStatement} of the core and paper-wide layers has a
theorem \texttt{X\_proved} whose type is literally \texttt{XStatement}
or, for a statement with parameters,
$\forall$\,\texttt{params, XStatement params}. A few details of the
corrected text are certified by auxiliary theorems of the formal proofs, which are kernel-checked but were not frozen in
advance. One of these details, the constant $20$ in Proposition~\ref{prop:inversion}, is also stated in three of the five
exceptions above; the other two state the symmetry-breaking sentence of Corollary~\ref{cor:insolv}. All five were frozen
after a proof of their content existed, and passed three independent reviews of their fidelity to the paper before their
\texttt{X\_proved} theorems were written. Two audit modules are part of the default build: one prints the
axiom set of every declaration it lists from the core layer and from
companion developments, and one discovers all frozen statements of the
paper-wide layer from the environment and fails the build unless each
has such a theorem depending only on the three standard axioms of the Lean/mathlib
foundation (propositional extensionality, choice, quotient soundness). The
development contains no unproved placeholder and no further axiom, and
every theorem and every safe definition in the 89 modules
reachable from the library's root module, the root included
(3036 theorems and 804 definitions), was additionally
re-checked by the Lean kernel.

Three features of the formalization deserve mention. First, it is
self-contained precisely where the prose leans on the literature: the
signed spanning-set count (Lemmas~\ref{lem:mobius}--\ref{lem:mobius2}), the Grassmannian
subspace count, and the $p$-group M\"obius evaluation (in its
subspace-lattice instance) are proved from first principles; so are the
ingredients of the second derivation of Section~\ref{sec:second-derivation}
--- Lemma~7(i) of \cite{cg2017} for every finite group, the product rule
for the M\"obius function of the subgroup lattice, and Hall's evaluation
on the primary components of $\A$ --- as well as the bounds of
\cite{gjs2016} and \cite{cg2017} specialized to $\ZZ_n^d$ and
\cite[Lemmas~E.21, E.22]{banda2019}. No verified step rests on the cited
results of Stanley, Hall, or Hawkes--Isaacs--\"Ozaydin. Second, the certified proof
of Theorem~\ref{thm:lower-bound} takes a shorter route than
Appendix~\ref{app:lower}: with the enumeration theorem available as a
verified identity, the prime case reduces to a reindexing of divisors of
$p^d$, and the composite case follows from the single fixed-set bound
$|S_{n^d}|\ge\sigma^{n^d/p_1}$ ($p_1$ the least prime factor) together with
elementary coefficient estimates (chiefly $\gaussbinom dr_q\le 2^d
q^{r(d-r)}$). In particular the certified chain nowhere invokes
\cite[Lemma~E.21]{banda2019} or its erratum, so the verified results are
independent of the published 2D bounds and their repairs. Third, and more
generally, the prose proofs are not replayed line by line: the formal
proofs follow their own routes, and proof-internal claims of the prose
proofs are formalized only selectively.

The Lean sources reside in \texttt{research/formal/shiftsym/} in the
project repository (the paper-wide layer in
\texttt{Shiftsym/Paper/} and \texttt{Shiftsym/AuditPaper.lean}), together
with the formalization plan and the
statement-preservation table.

\section{Summary and conclusion}\label{sec:conclusion}

\begin{table}[ht]
\centering
\caption{2D (\cite{banda2019}) vs.\ the $d$-dimensional generalization.}
\begin{tabular}{p{4.2cm}|p{4.2cm}|p{5.6cm}}
\toprule
\textbf{Component} & \textbf{2D} & \textbf{$d$-dimensional (this paper)} \\
\midrule
Lattice, cells & $\ZZ_n^2$, $N=n^2$ & $\ZZ_n^d$, $N=n^d$ \\
Order-$p$ subgroups per prime & $p+1$ & $(p^d-1)/(p-1)$ \\
Coefficients & $1,\;p{+}1,\;p$ & $r_d(p,v)=p^{\binom v2}\gaussbinom dv_p$ \\
Enumeration & Thm.~E.19 & Theorem~\ref{thm:main-enum}; finite abelian
groups: Theorem~\ref{thm:abelian} \\
$k$-active counts & Thm.~E.32 & Theorem~\ref{thm:kactive} \\
Lower bound (tight iff prime) & $(n{+}1)\sigma^n{-}n\sigma$ &
$F_d(n,\sigma)$, Theorem~\ref{thm:lower-bound} \\
Upper bound & $6\log_2(n)\sigma^{n^2/2}$ &
$(2^d{-}1{+}2\log_2 n)\sigma^{n^d/2}$ \\
Irreversibility & Thm.~E.5 & Theorem~\ref{thm:irrev} (equivariance
explicit) + Proposition~\ref{prop:rank} \\
Detection, worst case & $O(n^3)$ & $O(\mathrm{sopf}_d(n)\cdot N)$ \\
Detection, average case & $O(n^2)$; polylog claimed for random $n$ &
$O(\sum_{p\mid n}p^d)$; this bound is polylog (fixed $d$) only on a set
of $n$ of density $0$ (e.g.\ polylog-smooth $n$) \\
Decay rate of $P=\sigma^{-cN}$ & --- & $c\to1-1/p_1(n)$
(Proposition~\ref{prop:rate}) \\
Prime local minima of $P$ & yes for $p\ge5$ (empirical there; a theorem
here, Observation~\ref{obs:inversion}) & \textbf{fails for }$d\ge3$
(at $p=5$, $\sigma=2$); inversion of $(p,p{+}1)$ around
$d_c(p)\approx p\ln2$
(Proposition~\ref{prop:inversion}, Observation~\ref{obs:extrema}) \\
\bottomrule
\end{tabular}
\end{table}

The structure of shift-symmetry on $\ZZ_n^d$ is governed by the prime
factorization of $n$ through the subspace lattices of $\FF_p^d$: the
Gaussian binomial coefficients replace the ad hoc 2D coefficients, the
M\"obius values $(-1)^rp^{\binom r2}$ supply the signs, and every
dynamical statement flows from translation equivariance alone.
The $p$-torsion carries all of this: a configuration is shift-symmetric iff
some subgroup of prime order fixes it (Lemma~\ref{lem:completeness}), the
detector tests only these subgroups (Section~\ref{sec:detection}),
translation-equivariant dynamics never lower $\dim_{\FF_p}(\Stab(\ssb)\cap\A[p])$
(Proposition~\ref{prop:rank}), and on a finite abelian group the count depends only on
its order and the dimensions of its $p$-torsion (Corollary~\ref{cor:abelian-ranks}). Two
cautionary findings of the verification deserve
emphasis: naive dimensional lifting of true 2D statements can fail
outright (Remark~\ref{rem:oldbound-false},
Observation~\ref{obs:extrema}), and the enumeration principle itself,
though not the multi-prime closed form, was already available in the general
group-theoretic literature \cite{cg2017} --- both are reminders that the
value added here lies in the explicit formulas, the tightness
characterization, and the cellular-automata consequences, each verified
computationally at the definition level.

\paragraph{Applications.}
The counts use only a finite alphabet and a torus, and the dynamical
results use only translation equivariance (their proofs do not need a finite
alphabet), so they apply wherever these ingredients occur.
\begin{itemize}[leftmargin=*]
    \item \emph{Distributed computing.} A CA on $\ZZ_n^d$ is a synchronous anonymous
    network (a torus), a special case of the networks for which Angluin proved
    the classical impossibility of symmetry breaking \cite{angluin1980};
    Yamashita and Kameda characterized when leader election and related
    problems are solvable on anonymous networks \cite{yamashita1996}.
    For CA, the symmetry obstruction and its use in performance bounds were
    studied in one dimension by Banda, Caughman, and Posp\'ichal
    \cite{banda2015}, and in two dimensions by Banda et al.~\cite{banda2019}.
    Corollary~\ref{cor:insolv} and the counts extend this obstruction and its
    quantitative form to arbitrary dimension: for $N>1$, every translation-equivariant CA
    fails at leader election on all inputs in $S_{n^d}$, a fraction
    $P^{\mathrm{unif}}_{n^d}$ of all inputs, bounded in
    Corollary~\ref{cor:prob}.
    \item \emph{Machine learning.} A stride-one convolutional layer with circular
    padding, shared weights, and pointwise nonlinearities is a
    translation-equivariant map of feature maps on $\ZZ_n^d$ \cite{cohen2016}.
    Since the proof of Theorem~\ref{thm:irrev} uses equivariance alone, no stack
    of such layers can remove a translation symmetry of its input, an instance
    of the limitation of equivariant networks discussed by Kaba and Ravanbakhsh
    \cite{kaba2023}. Neural cellular automata \cite{mordvintsev2020} grow
    patterns from a single seed cell, which is an asymmetric configuration, and
    update their cells at random; in the terms of Corollary~\ref{cor:insolv},
    they restrict the inputs and drop determinism and synchrony.
    \item \emph{Materials science.} Enumerations of the derivative superstructures of
    an alloy discard the ``super-periodic'' labelings, whose period is shorter
    than the supercell \cite{hartforcade2008}. For a Bravais parent lattice,
    $\sigma$ atomic species, and the supercell spanned by $n$ times the parent
    basis vectors, these are exactly the configurations of $S_{n^d}$ ($d=3$ for
    bulk crystals), so Theorem~\ref{thm:main-enum} counts them, and for binary
    alloys Theorem~\ref{thm:kactive} counts them at fixed composition. A general
    supercell has a quotient of the form $\ZZ_{a}\times\ZZ_{b}\times\ZZ_{c}$
    (Smith normal form), a rectangular torus, whose super-periodic labelings
    Theorem~\ref{thm:abelian} counts (for binary alloys at fixed
    composition, Remark~\ref{rem:abelian}(i)).
    \item \emph{Physics.} In exact diagonalization of quantum lattice models with
    translation symmetry, basis states whose period is shorter than the lattice
    need a modified normalization and are compatible with only some momenta
    \cite[\S\S4.1.3, 4.4.1]{sandvik2010}; on $\ZZ_n^d$ with $\sigma$ local states, these
    are exactly the configurations of $S_{n^d}$. More generally, by the argument
    of Theorem~\ref{thm:irrev}, any deterministic update map on a periodic
    grid that commutes with the grid translations (for example, a synchronous
    update applying the same rule at every site, as in a lattice gas or an
    explicit discretization of a reaction--diffusion equation) preserves, in
    exact arithmetic, every translation symmetry of its initial data; symmetric
    initial data must therefore be perturbed (round-off can act as such a
    perturbation) for symmetry breaking to be observed in such simulations.
    \item \emph{Random initial states.} In cellular-automaton pseudorandom generators
    \cite{wolfram1986} and in simulations started from random configurations, a
    seed in $S_{n^d}$ confines the whole trajectory to $S_{n^d}$
    (Theorem~\ref{thm:irrev}). For a uniform seed this happens with probability
    at most $(2^d-1+2\log_2 n)\,\sigma^{-n^d/2}$ (Corollary~\ref{cor:prob}),
    which is negligible except on very small lattices. Structured or sparse
    seeds are another matter (the all-zero seed lies in $S_{n^d}$). For a
    uniform seed with a prescribed number $k$ of active cells,
    Theorem~\ref{thm:kactive} gives the exact probability
    $|S^a_{n^d,k}|/\bigl(\binom Nk(\sigma-1)^{N-k}\bigr)$, and
    Corollary~\ref{cor:dens} gives it when $k$ is itself uniform; the
    algorithm of Section~\ref{sec:detection} tests a seed in
    $O(\mathrm{sopf}_d(n)\cdot n^d)$ cell reads.
\end{itemize}

\paragraph{Verification and its audit.}
A formal proof certifies a paper only as far as its formal statements say
what the paper says: the kernel checks the proofs, not the reading. Besides
the fidelity reviews of Section~\ref{sec:formal}, the correspondence between
every theorem, lemma, proposition, and corollary of this paper and its
frozen Lean statements was therefore audited by hand by one of us (T.L.),
with a viewer built for the purpose.\footnote{The viewer and the audit log are in \texttt{research/formal/correspondence\_audit/} of the project repository;
T.L. is developing a product version of the viewer (\url{https://functorial.inc/}).} The viewer shows the official PDF beside the
Lean declarations and links sub-expressions of the paper to the Lean
sub-expressions claimed to express them ($529$ links for the $27$ items),
each link classified by kind: identical, a change of representation, a
definition, a hypothesis, or a generalization (a few links point to a
citation or have a counterpart on one side only). No link is of the kind
``weaker'', which would mark a part of the paper not covered by a frozen
statement. The highlighted regions are computed from an instrumented copy
of the paper whose words must typeset at exactly the positions they have in
the official PDF; every action of the auditor is recorded, with timestamps,
in a log kept under version control; and when the viewer data are rebuilt,
a verdict is flagged as stale if the item's \TeX{} source, the Lean source of
its declarations, or its link data differ from those it was given against.
All $27$ items were approved. A verdict given on an earlier draft was carried
over by a script only where the item's \TeX{} source, the Lean source of its
declarations, and its links were unchanged apart from line numbers; otherwise
the item was judged again.
A separate check of the link data, made during the audit, found
that one clause, the symmetry-breaking sentence of
Corollary~\ref{cor:insolv}, which the links had classified as having no Lean
counterpart, was not covered by any frozen statement; it was reworded and
frozen, and its new links were then audited. Tools of this kind make the
alignment between informal and formal statements inspectable item by item,
and we regard them as a natural companion to proof assistants for papers
like this one.

\paragraph{Open problems.}
(1) A theorem (not only computation) delimiting for which $(n,d)$ the
prime-dip/even-peak pattern of $P^{\mathrm{unif}}$ holds. (2) Rectangular
tori $\ZZ_{n_1}\times\dots\times\ZZ_{n_d}$, and finite abelian groups in
general (counted by Theorem~\ref{thm:abelian}): a lower envelope for the
number of shift-symmetric configurations that reduces to $F_d(n,\sigma)$ of
Theorem~\ref{thm:lower-bound} on $\ZZ_n^d$ and is attained exactly on a
natural class of groups. (For a fixed order alone, the minimum is attained by the cyclic
group, Corollary~\ref{cor:abelian-ranks}.)

\section*{Acknowledgments}
AI assistants, including but not limited to Claude (Anthropic) and Codex (OpenAI), were used in developing and reviewing the proofs, the Lean formalization, and the verification scripts, and in editing the text; the authors take full responsibility for the content.

\appendix\section*{Appendices}\addcontentsline{toc}{section}{Appendices}
\section{Proof of Theorem~\ref{thm:lower-bound}}\label{app:lower}

Throughout, $q^m(c)$ denotes the truncation of \eqref{eq:main-formula} to
the first $m$ primes with exponents scaled by $c$:
\[
  q^m(c)\;=\;\sum_{\substack{\vv\in\{0,\dots,d\}^m\\\vv\neq0}}
  (-1)^{1+|\vv|}\prod_{i=1}^m r_d(p_i,v_i)\,
  \sigma^{\,cN\prod_i p_i^{-v_i}},\qquad q^m:=q^m(1),
\]
so $q^k=|S_{n^d}|$ by Theorem~\ref{thm:main-enum}. We write $p:=p_{m+1}$,
$N_j:=n^d/p^j$, $t:=n/p$, $D:=d(d-1)/2$.

If $n$ is prime, Theorem~\ref{thm:lower-bound} holds with equality by
Corollary~\ref{cor:prime}. For composite $n$ we prove
$q^1\le q^2\le\dots\le q^k$ (Lemma~\ref{lem:chain}) and
$q^1>F_d(n,\sigma)$ (Lemma~\ref{lem:base} for $d\ge3$; $d\le2$ below),
which give strict inequality at every composite $n$.

\begin{lemma}[Union interpretation on rectangular tori]\label{lem:dtorus}
Let $1\le m<k$ and $1\le j\le d$. Then
\[
  0\;\le\;q^m(p^{-j})\;\le\;\sigma^{N_j},
  \qquad\text{and}\qquad
  q^m(p^{-1})\;\le\;\sum_{i\le m}\gaussbinom d1_{p_i}\,\sigma^{N_1/p_i}.
\]
\end{lemma}

\begin{proof}
Consider $T_j=(\ZZ/n\ZZ)^{d-j}\times(\ZZ/\tfrac np\ZZ)^{j}$ with
$N_j$ cells. For $i\le m$, $p_i$ divides both $n$ and $n/p$ (as
$p_i\neq p$), so the $p_i$-torsion of $T_j$ is again $\cong\FF_{p_i}^d$
and Lemmas \ref{lem:fix}, \ref{lem:torsion}--\ref{lem:joins}, \ref{lem:mobius} and \ref{lem:mobius2} apply
verbatim on $T_j$. By the computation in the proof of
Theorem~\ref{thm:main-enum}, restricted to $p_1,\dots,p_m$,
$q^m(p^{-j})$ equals the number of configurations of $T_j$ fixed by some
order-$p_i$ cyclic subgroup, $i\le m$ --- a union of fixed-point sets.
Nonnegativity, the trivial bound $\sigma^{N_j}$, and Boole's inequality
(at $j=1$) follow. Exponent integrality:
$v_{p_i}(N_j)=d\,\alpha_{p_i}\ge d\ge v_i$.
\end{proof}

\begin{lemma}[Coefficient estimates]\label{lem:coef}
Let $q\ge2$ be real, with $\gaussbinom dj_q$ and
$r_d(q,j)=q^{\binom j2}\gaussbinom dj_q$ defined by the product formula.
Then $\gaussbinom dj_q\le 4q^{\,j(d-j)}$,
$\gaussbinom d1_q\le 2q^{\,d-1}$, and
$\sum_{j=2}^{d}r_d(q,j)\le 4(d-1)\,q^{D}$.
\end{lemma}

\begin{proof}
$\gaussbinom dj_q=q^{j(d-j)}\prod_{i=1}^{j}
\frac{1-q^{-(d-j+i)}}{1-q^{-i}}\le
q^{j(d-j)}\prod_{i\ge1}(1-2^{-i})^{-1}$, and
$\prod_{i\ge1}(1-2^{-i})\ge\frac12\cdot\frac34\cdot\frac78
\bigl(1-\sum_{i\ge4}2^{-i}\bigr)=\frac{147}{512}>\frac14$; also
$\gaussbinom d1_q\le q^d/(q-1)\le2q^{d-1}$. Hence
$r_d(q,j)\le4q^{\,j(2d-j-1)/2}$; the exponent is a concave parabola in
$j$ with vertex $d-\tfrac12$, so over $2\le j\le d$ its maximum is $D$,
attained at $j=d$ (and $j=d-1$ for $d\ge3$); the tail sum is at most
$4(d-1)q^D$.
\end{proof}

\begin{lemma}[Monotone chain]\label{lem:chain}
For composite $n$ with $k\ge2$, $d\ge1$, $\sigma\ge2$, and $1\le m<k$:
$q^m\le q^{m+1}$.
\end{lemma}

\begin{proof}
Splitting $q^{m+1}$ by the value $j$ of $v_{m+1}$ and isolating, for
$j\ge1$, the term with $(v_1,\dots,v_m)=\mathbf 0$:
\begin{equation}\label{eq:decomp}
  q^{m+1}-q^m\;=\;\sum_{j=1}^{d}(-1)^{j+1}r_d(p,j)\,B_j,
  \qquad B_j:=\sigma^{N_j}-q^m(p^{-j})\;\in\;[0,\sigma^{N_j}]
\end{equation}
by Lemma~\ref{lem:dtorus}. For $d=1$ the sum is the single term
$B_1\ge0$; assume $d\ge2$. Since $m\ge1$: $p=p_{m+1}\ge3$, and
$p\prod_{i\le m}p_i\mid n$ gives $p\le n/2$ and $t\ge2$. We use:
\begin{itemize}
\item[(I)] $\sigma^{N_1/2}\ge 4\,t^{\,d-1}$;
\item[(II)] $\sigma^{N_1-N_2}\ge 8(d-1)\,p^{D}$.
\end{itemize}
Given (I): by Lemma~\ref{lem:dtorus} and
$\sum_{i\le m}\gaussbinom d1_{p_i}\le\sum_i2p_i^{d-1}
\le2\bigl(\sum_ip_i\bigr)^{d-1}\le2\bigl(\prod_ip_i\bigr)^{d-1}
\le2t^{\,d-1}$ (superadditivity of $x^{d-1}$; $\sum p_i\le\prod p_i$ for
distinct primes; $\prod_{i\le m}p_i\le n/p=t$),
\[
  B_1\;\ge\;\sigma^{N_1}-2t^{\,d-1}\sigma^{N_1/2}\;\ge\;
  \tfrac12\sigma^{N_1}.
\]
Given (II): dropping the positive $j\ge3$ odd terms of \eqref{eq:decomp},
bounding $B_j\le\sigma^{N_2}$ for $j\ge2$, and using
Lemma~\ref{lem:coef} with $r_d(p,1)\ge1$:
$q^{m+1}-q^m\ge\tfrac12\sigma^{N_1}-4(d-1)p^{D}\sigma^{N_2}\ge0$.

It remains to prove (I), (II) for composite $n\ge6$, $d\ge2$. (I):
$N_1/2=\tfrac t2n^{d-1}\ge n^{d-1}$ and $t\le n/3$, so it suffices that
$n^{d-1}\ge2+(d-1)\log_2(n/3)$; since $\log_2 n\le n/2$ for $n\ge4$, it
suffices that $n^{d-1}\ge2+(d-1)n/2$ --- for $d=2$ this is $n\ge4$; for
$d\ge3$, $n^{d-1}\ge6^{d-2}n\ge(\tfrac{d-1}2)n+2$ since
$6^{d-2}\ge\tfrac{d-1}2+1$ (at $d=3$: $6\ge2$; per increment of $d$ the
left side sextuples while the right grows by $\tfrac12$). (II):
$N_1-N_2=N_1(1-\tfrac1p)\ge\tfrac23N_1\ge\tfrac43n^{d-1}$ and
$p^D\le(n/2)^D$, so it suffices that
$\tfrac43n^{d-1}\ge3+\log_2(d-1)+D\log_2(n/2)$. For $d=2$:
$\tfrac43n\ge3+\log_2(n/2)$ holds for $n\ge6$. For $d\ge3$, using
$\log_2(d-1)\le d-1$ and $\log_2(n/2)\le n/2$, it suffices that
$\tfrac43\,6^{d-2}n\ge2+d+\tfrac{d(d-1)}4n$, which holds for $n\ge6$
because $\tfrac43\,6^{d-2}\ge\tfrac{d(d-1)}4+\tfrac{2+d}6$: at $d=3$
this reads $8\ge\tfrac32+\tfrac56$, and the induction step
$6\bigl(\tfrac{d(d-1)}4+\tfrac{2+d}6\bigr)\ge
\tfrac{(d+1)d}4+\tfrac{3+d}6$ reduces to $5d\ge7$.
\end{proof}

\begin{lemma}[Base case, $d\ge3$]\label{lem:base}
For composite $n$, $\sigma\ge2$, $d\ge3$: $q^1>F_d(n,\sigma)$.
\end{lemma}

\begin{proof}
Let $p_1$ be the smallest prime factor; $p_1\le n/2$ and $p_1^2\le n$
(for $n=p_1^{\alpha}$, $\alpha\ge2$; otherwise $p_1^2<p_1q\le n$ for the
second prime $q$). Dropping positive $r\ge3$ terms of $q^1$ and negative
terms of $F_d$, and using Lemma~\ref{lem:coef} at $q=p_1$ and at the
composite argument $q=n$ together with
$\sigma^{n^d/p_1^r}\le\sigma^{n^d/p_1^2}$ ($r\ge2$) and
$\sigma^{n^{d-r}}\le\sigma^{n^{d-3}}$ ($r\ge3$ odd):
\[
  q^1\ge\gaussbinom d1_{p_1}\sigma^{n^d/p_1}
  -4(d-1)p_1^{D}\sigma^{n^d/p_1^2},
  \qquad
  F_d\le\gaussbinom d1_n\sigma^{n^{d-1}}
  +4(d-1)n^{D}\sigma^{n^{d-3}} .
\]
All three exponent gaps are at least $n^{d-1}$:
$\tfrac{n^d}{p_1}-\tfrac{n^d}{p_1^2}\ge\tfrac{n^d}{2p_1}\ge n^{d-1}$ and
$\tfrac{n^d}{p_1}-n^{d-1}\ge n^{d-1}$ (as $n/p_1\ge2$), a fortiori for
$n^{d-3}$. With $\gaussbinom d1_{p_1}\ge3$, $\gaussbinom d1_n\le
2n^{d-1}$, $p_1\le n$:
\[
  q^1-F_d\;\ge\;\sigma^{n^d/p_1}\Bigl[3-\bigl(8(d-1)n^{D}
  +2n^{d-1}\bigr)2^{-n^{d-1}}\Bigr]\;\ge\;\sigma^{n^d/p_1}>0,
\]
provided $2^{n^{d-1}}\ge12d\,n^{D}$ (using
$8(d-1)n^D+2n^{d-1}\le12dn^D$, as $n^{d-1}\le n^D$ for $d\ge2$). Taking
logarithms it suffices that $n^{d-1}\ge4+d+\tfrac{d(d-1)}2\log_2n$ for
$n\ge4$, $d\ge3$: at $(4,3)$ it reads $16\ge13$; for fixed $d$ the margin
increases in $n$ (derivative $(d-1)n^{d-2}-\tfrac{d(d-1)}{2n\ln2}>0$),
and the step $d\to d{+}1$ adds $1+d\log_2n\le1+n^{d-1}\le(n-1)n^{d-1}$
to the right side while multiplying the left by $n$.
\end{proof}

\begin{proof}[Proof of Theorem~\ref{thm:lower-bound}]
$d=1$: $F_1=\sigma$; primes give $|S_n|=\sigma$, composites
$|S_n|\ge\sigma^{n/p_1}\ge\sigma^2>\sigma$.
$d=2$: $F_2(n,\sigma)=(n{+}1)\sigma^n-n\sigma$ and the claim is
\cite[Lemma~E.21]{banda2019}; we import its base step
($q^1>(n{+}1)\sigma^n-n\sigma$ at composite $n$, final display of the
source proof; the inequality holds at every composite $n$, including the
prime powers $p^\alpha$, $\alpha\ge2$, which the case split of the source
(``If $k=1$, i.e., $n$ is a prime'') does not treat (corrected in the companion erratum), and it is
machine-checked, Section~\ref{sec:formal}), while the chain step is re-proved uniformly by
Lemma~\ref{lem:chain} with $d=2$ (the published chain, Lemma~A.21 of
\cite{banda2019}, has three exceptional lattice sizes, repaired in the
companion erratum).
$d\ge3$: prime $n$ gives equality (Corollary~\ref{cor:prime}); composite
$n$ with $k=1$ gives $|S|=q^1>F_d$ (Lemma~\ref{lem:base}); composite $n$
with $k\ge2$ gives $|S|=q^k\ge\dots\ge q^1>F_d$
(Lemmas~\ref{lem:chain}, \ref{lem:base}).
\end{proof}

Every inequality in this appendix was verified exactly on parameter grids
before submission (decomposition identity, torus sandwich, union bound,
coefficient estimates, conditions (I)/(II), the chain, and the base case;
see the reproducibility paragraph of Section~\ref{sec:empirical}), and the
argument passed independent adversarial review.

\section{Finite abelian groups: dependence on the ranks}\label{app:abelian}

Notation is as in Section~\ref{sec:abelian}.

\begin{corollary}[Dependence on the ranks]\label{cor:abelian-ranks}
Let $\sigma\ge2$, and write $\mathrm{ac}(G)=\sigma^{|G|}-|S_G|$ for the number of
aperiodic configurations.
\begin{enumerate}[label=(\roman*)]
\item $|S_G|$ depends on $G$ only through $N$ and its $p$-ranks; for
example, $\ZZ_2\times\ZZ_{24}$ and $\ZZ_4\times\ZZ_{12}$ both have $50\,388\,256$
shift-symmetric configurations for $\sigma=2$.
\item If $G'$ has the same order as $G$, $d_p(G')=d_p(G)+1$ for one prime
$p$, and the other ranks agree, then
$|S_{G'}|-|S_G|=p^{d_p(G)}\,\mathrm{ac}(G'')>0$ for any abelian group $G''$ of
order $N/p$ with the ranks of $G$. In particular, among the abelian groups
of order $N$, the cyclic group has the fewest shift-symmetric
configurations.
\end{enumerate}
\end{corollary}

\begin{proof}
(i) is immediate from Theorem~\ref{thm:abelian}. (ii) In complement form,
\[
  \mathrm{ac}(G)=\sum_{\mathbf 0\le\vv\le(d_p)_p}\ \prod_p c(p,d_p,v_p)\,
  \sigma^{N/\prod_p p^{v_p}},\qquad
  c(p,e,v)=(-1)^vp^{\binom v2}\gaussbinom ev_p .
\]
The $q$-Pascal rule
$\gaussbinom{e+1}v_p=\gaussbinom ev_p+p^{e+1-v}\gaussbinom e{v-1}_p$ gives
\[
  c(p,e,v)-c(p,e+1,v)=p^e\,c(p,e,v-1),\qquad c(p,e,e+1)=0 .
\]
Subtracting the two sums term by term and shifting $v_p$ by one gives
\[\mathrm{ac}(G)-\mathrm{ac}(G')=p^{d_p(G)}\,\mathrm{ac}(G'').\]
Finally, $\mathrm{ac}(G'')\ge1$: the indicator of a single cell has trivial
stabilizer when $|G''|\ge2$, and $\mathrm{ac}(G'')=\sigma$ when $|G''|=1$. Write
$\alpha_p$ for the exponent of $p$ in $N$. An abelian group of order $N$ has
$1\le d_p\le\alpha_p$ for every $p\mid N$, and for $1\le r\le\alpha$ the group
$\ZZ_{p^{\alpha-r+1}}\times\ZZ_p^{r-1}$ has order $p^\alpha$ and $p$-rank $r$; so
$G''$ exists (as $d_p(G)+1\le\alpha_p$), every rank profile with
$1\le d_p\le\alpha_p$ occurs, and every abelian group of order $N$ is
reached from the cyclic group, which has $d_p=1$ for every $p\mid N$, by
such steps.
\end{proof}

\begin{remark}\label{rem:abelian}
(i) The $k$-active counts of Theorem~\ref{thm:kactive} extend to $G$ in the
same way, with $r_d$ replaced by $r_{d_p}$ and $n^d$ by $N$.
(ii) $G$ acts freely on the aperiodic configurations, so
$N\mid\sigma^N-|S_G|$, and $(\sigma^N-|S_G|)/N$ is the number of
$G$-orbits of aperiodic configurations; for rectangular tori these are the
atranslational necklaces of \cite{adamson2021}.
\end{remark}


\end{document}